\documentclass[2p]{article}
\usepackage{amssymb,amsmath,amsthm}
\usepackage{indentfirst}
\usepackage{exscale}
\usepackage{relsize}
\usepackage{xcolor}
\usepackage[numbers,sort&compress]{natbib}

\usepackage{geometry}

\theoremstyle{definition}

\allowdisplaybreaks
\theoremstyle{remark}

\numberwithin{equation}{section}

\begin{document}

\title{\Large\bf{Two nontrivial solutions for a nonhomogeneous quasilinear elliptic system  with sign-changing weight functions}
 }
\date{}
\author {\ Wanting Qi$^{1}$, \ Xingyong Zhang$^{1,2}$\footnote{Corresponding author, E-mail address: zhangxingyong1@163.com}\\
      {\footnotesize $^{1}$Faculty of Science, Kunming University of Science and Technology, Kunming, Yunnan, 650500, P.R. China.}\\
      {\footnotesize $^{1,2}$Research Center for Mathematics and Interdisciplinary Sciences, Kunming University of Science and Technology,}\\
 {\footnotesize Kunming, Yunnan, 650500, P.R. China.}\\
 }
 \date{}
 \maketitle

 \begin{center}
 \begin{minipage}{15cm}

\small  {\bf Abstract:}
We are interested in looking for two nontrivial solutions for a class of nonhomogeneous quasilinear elliptic system with sign-changing weight functions and concave-convex nonlinearities on the bounded domain.
This kind of quasilinear elliptic system arises from nonlinear optics, whose feature is that its differential operator depends on not only $\nabla u$  but also $u$.
Employing the mountain pass theorem and Ekeland's variational principle as the major tools,
we show that the system has at least one nontrivial solution of positive energy and one nontrivial solution of negative energy, respectively.
\end{minipage}
 \end{center}
  \allowdisplaybreaks
 \vskip2mm
\section{Introduction and main results}
\par
Mathematically, the guided TM-modes, which propagating through a self-focusing anisotropic dielectric, are special solutions of Maxwell's equations with a nonlinear constitutive relation of a type commonly used in nonlinear optics when treating the propagation of waves in a cylindrical wave-guide.
Under a series of transformations, the existence of such modes propagating in a nonlinear dielectric medium can be reduced to searching for solutions of the corresponding one-dimensional nonlinear eigenvalue problem on the Sobolev space $H_{0}^{1}(0,\infty)$.
For details, we infer the readers to \cite{Landau1984,Mihalache1989,Svelto1974} for the essential physical background, and to \cite{Chen1991,Snyder1991,Stuart1996,Stuart1997,Stuart2010,Stuart2001} for the standard procedure to study TM-modes in nonlinear optics.
Specifically, in \cite{Stuart2001}, Stuart-Zhou considered the following nonlinear eigenvalue problem:
\begin{equation}\label{2001.1.0}
\left\{
  \begin{array}{ll}
-\left\{\gamma\left(\frac{1}{2}\left[w^{2}+\left(w'+\frac{w}{r}\right)^{2}\right]\right)\left(w'+\frac{w}{r}\right)\right\}'
+\gamma\left(\frac{1}{2}\left[w^{2}+\left(w'+\frac{w}{r}\right)^{2}\right]\right)w
=\lambda w,\\
w(0)=0, \lambda=(\frac{\omega}{kc})^{2},
  \end{array}
 \right.
 \end{equation}
where $ r$ is the distance from the axis, $c>0$ is the speed of light in a vacuum, $k$ is a real number related to the wavelength, $\omega$ is the frequency of the magnetic field, and $ \gamma:[0,+\infty)\rightarrow(0,+\infty)$ satisfies some reasonable conditions.
Using the mountain pass theorem, they established the existence of nontrivial solutions for problem (\ref{2001.1.0}), and showed that these solutions can generate solutions of Maxwell's equations which having the form of guided travelling waves propagating through a homogeneous self-focusing dielectric.
\par
Inspired by these works described earlier, in \cite{Stuart2011}, Stuart introduced the following second-order quasilinear elliptic problem:
\begin{equation}\label{2011.1.0}
 \left\{
  \begin{array}{ll}
  -\mbox{div}\left\{\phi \left(\frac{u^{2}+|\nabla u|^{2}}{2}\right)\nabla u\right\}+\phi\left (\frac{u^{2}+|\nabla u|^{2}}{2}\right)u
  =
  \lambda u+h, &\mbox{in } \Omega,\\
 u=0, &\mbox{on } \partial\Omega,
  \end{array}
 \right.
 \end{equation}
which has some similar but more simple features to problem (\ref{2001.1.0}),
where $\Omega$ is a bounded open subset of $\mathbb{R}^{N}$, $N\geq1$, $\lambda\in\mathbb{R}$, $\phi:[0,+\infty)\rightarrow\mathbb{R}$ is a positive continuous function which satisfies the monotonicity condition and other appropriate conditions.
Under the assumption that $h\in L^{2}(\Omega)$ and $h\geq0$ a.e. on $\Omega$,
they obtained the existence of two non-negative weak solutions for problem (\ref{2001.1.0}) : one is the local minimum of the corresponding functional and the other is mountain-pass solution.
Afterwards, in \cite{Jeanjean2022}, Jeanjean-R\u{a}dulescu considered the existence of solutions for equation (\ref{2011.1.0}) with nonlinearities $f(u)+h$, which extend $\lambda u+h$ directly to the case of the nonlinear growth reaction term.
To be specific, $h$ is non-negative and $f$ is a given continuous function which has either a sublinear decay or  a linear growth at infinity.
In the sublinear decay case, $\phi$ was assumed to satisfy some reasonable conditions, and the existence of non-negative solutions for problem (\ref{2011.1.0}) was established by making use of a minimization procedure.
In the linear growth case, by means of the mountain pass theorem,
they proved problem (\ref{2011.1.0}) has at least one or two non-negative solutions.
In \cite{Pomponio2021}, Pomponio-Watanabe studied problem (\ref{2011.1.0}) with general nonlinear terms of Berestycki-Lions' type in the entire space $\mathbb{R}^{N}$.
They proved that,
when the nonlinearity is a continuous and odd function, which has at least a linear growth near the origin,
problem (\ref{2011.1.0}) possesses at least a non-trivial weak solution.
The proof relies upon the mountain pass theorem and a technique of adding one dimension for space $\mathbb{R}^{N}$.
Then by establishing the regularity of solutions and the Pohozaev identity, they also obtained the existence of a radial ground state solution and a ground state solution.
\par
It is easy to see that the assumptions on reaction term in \cite{Jeanjean2022,Pomponio2021} prohibit the consideration of a large class of functions which behaves sublinearly near the origin and superlinearly near infinity .
For example, $\lambda |u|^{q-2}u+|u|^{p-2}u$ with $\lambda>0$ and $1<q<2<p\leq 2^{\ast}(2^{\ast}=2N/N-2 \;\mbox{if}\; N>2; 2^{\ast}=\infty \;\mbox{if}\; N=1,2)$.
This kind of nonlinearity is usually called the type of concave-convex, which has been initially treated by Ambrosetti-Brezis-Cerami \cite{Ambrosetti1994}.
In \cite{Ambrosetti1994}, by means of the subsupersolutions and variational method, existence, multiplicity, and nonexistence of positive solutions was obtained for the semilinear elliptic problem:
\begin{equation}\label{1994.1.0}
 \left\{
  \begin{array}{ll}
 -\Delta u=\lambda a(x)|u|^{q-2}u+b(x)|u|^{p-2}u, &\mbox{in } \Omega,\\
 u=0, &\mbox{on } \partial\Omega,
  \end{array}
 \right.
 \end{equation}
with $a\equiv b\equiv1$ and $u>0$ in $\Omega$, where $\Omega$ is a bounded regular domain of $\mathbb{R}^{N}$($N>2$) and $\lambda$ is a positive parameter.
For the relationship of the number of solutions and the sign of the weight functions: $a,b:\Omega\rightarrow\mathbb{R}$, further studies were taken by Wu \cite{Wu2006,Wu2009} and Brown-Wu \cite{Brown2007}.
For details, when the weight functions $b\equiv1$ and $a$ is a continuous function which change sign in $\Omega$,
the existence of two positive solutions for problem (\ref{1994.1.0}) were established in \cite{Wu2006}.
When the weight functions $a$ and $b$ are both change sign in $\Omega$, the multiplicity results of solutions for problem (\ref{1994.1.0}) were also obtained in \cite{Wu2009} and \cite{Brown2007}, respectively.
After the above mentioned works, the solvability of problem like (\ref{1994.1.0}) with sign-changing weight functions has been extensively studied, for example,
\cite{Wu2009,Chen2013,Wu2010,Chen2011,Carvalho2017,Silva2023,Afrouzi2009,Chen2022,Fan2013,Sahu1991,Brown2009} and the references therein.
Particularly, in \cite{Chen2013}, Chen-Huang-Liu investigated the multiplicity of solutions for the $p$-Kirchhoff elliptic problem with concave-convex and nonhomogeneous term on unbounded domain.
By applying the Mountain pass theorem and Ekeland's variational principle,
they proved that problem has at least one positive energy solution and one negative energy solution.
Subsequently, the idea in \cite{Chen2013} has been applied to a lot of different problems, for examples,
the fractional $p$-Kirchhoff problems \cite{Xiang2015,Zuo2019},
the Kirchhoff type problem involving a nonlocal operator \cite{Chen2016},
$p$-biharmonic equation \cite{Xiu2016},
sixth-order nonhomogeneous $p(x)$-Kirchhoff problems \cite{Hamdani2021} and so on.
\par
Recently, in \cite{Qi2023}, we extended the methods of \cite{Wu2009,Brown2007} to consider the multiplicity of solutions for problem (\ref{2011.1.0}) with concave-convex nonlinearities and sign-changing weight functions.
By the Nehari manifold and doing a fine analysis associated on the fibering map,
we obtained that problem (\ref{2011.1.0}) admits at least one positive energy solution  and negative energy solution  which is also the ground state solution of problem (\ref{2011.1.0}).
\par
Motivated by \cite{Chen2013,Qi2023}, in the present paper, we investigate the multiplicity results of solutions for the following quasilinear elliptic system with concave-convex and nonhomogeneous terms:
 \begin{equation}\label{eq1}
 \left\{
  \begin{array}{ll}
 -\mbox{div}\left\{\phi_{1} \left(\frac{u^{2}+|\nabla u|^{2}}{2}\right)\nabla u\right\}+\phi_{1}\left (\frac{u^{2}+|\nabla u|^{2}}{2}\right)u
 =
 \lambda a(x)|u|^{q-2}u+\frac{\alpha}{\alpha+\beta}b(x)|u|^{\alpha-2}u|v|^{\beta}+h_{1}(x), &\mbox{in } \Omega,
 \\
 -\mbox{div}\left\{\phi_{2} \left(\frac{v^{2}+|\nabla v|^{2}}{2}\right)\nabla v\right\}+\phi_{2}\left (\frac{v^{2}+|\nabla v|^{2}}{2}\right)v
 =
 \mu c(x)|v|^{q-2}v+\frac{\beta}{\alpha+\beta}b(x)|u|^{\alpha}|v|^{\beta-2}v+h_{2}(x), &\mbox{in } \Omega,\\
 u=v=0, &\mbox{on } \partial\Omega,
  \end{array}
 \right.
 \end{equation}
where $\Omega$ is a bounded domain in $\mathbb{R}^{N}( N\geq3)$ with smooth boundary,
$\lambda$ and $\mu$ are two positive parameters,
$a,b,c,h_{1},h_{2}:\Omega\rightarrow \mathbb{R}$ are continuous functions which may change sign on $\Omega$,
$\alpha,\beta>1$ and $1<q<2<\alpha+\beta<2^{\ast}$,
$\phi_{i}:[0,+\infty)\rightarrow \mathbb{R}$ $(i=1,2)$ are two continuous functions which satisfy the following conditions:
\begin{itemize}
 \item[$(\phi_1)$] there exist two constants $0< \rho_0\leq \rho_1$ such that
 $ 0<\rho_0\leq \phi_{i}(s)\leq \rho_1$ for all $s\in [0,+\infty)$;
 \item[$(\phi_2)$] the map $s\rightarrow \Phi_{i}(s^{2})$ is strictly convex on $[0,+\infty)$, where $\Phi_{i}(s):=\int_{0}^{s}\phi_{i}(\tau)d\tau$;
 \item[$(\phi_3)$] there exists a positive constant $\phi_{i}(\infty)$ such that $\phi_{i}(s)\rightarrow\phi_{i}(\infty)$ as $s\rightarrow +\infty$;
 \item[$(\phi_4)$] $\Phi_{i}(s)\geq \phi_{i}(s)s$ for all $s\in [0,+\infty)$.
 \end{itemize}
Moreover, we introduce the following assumptions on $a,b,c,h_{1}$ and $h_{2}$:
\begin{itemize}
 \item[$(A)$] there exist $\tau_{1}>0$ and $\tau_{2}>0$, such that
 $a\in L^{\frac{\alpha+\beta}{\alpha+\beta-q}+\tau_{1}}(\Omega)$ and
 $c \in L^{\frac{\alpha+\beta}{\alpha+\beta-q}+\tau_{2}}(\Omega)$;
 \item[$(B)$] $b\in L^{\infty}(\Omega)$ and there exists a non-empty open domain $\Omega_{1}\subset\Omega$ such that $b(x)>0$ in $\Omega_{1}$;
 \item[$(H)$] $h_{1},h_{2}\in L^{2}(\Omega)$ and $h_{1},h_{2}\not\equiv0$.
\end{itemize}
By the mountain pass theorem and Ekeland's variational principle, we obtain the existence of two nontrivial solutions for system (\ref{eq1}).
Our results develop those in some known references in the following sense:
\begin{itemize}
\item[(\uppercase\expandafter{\romannumeral1})]
At first, for the system like (\ref{eq1}), as far as we know, there is no paper to study the multiplicity of solutions by the variational method.
Moreover, when $\phi_1\equiv\phi_2\equiv1$ and $h_{1}\equiv h_{2}\equiv0$, the problem (\ref{eq1}) reduces to the semilinear elliptic system with concave-convex terms and sign-changing weight functions, which has been studied extensively over the past several decades (for example, see \cite{Wu2008,Lin2012,Fan2014,Fan2015,Chen2012,Li2013,Li2014} and references therein).
However, there seems to be little work for the case $h_{1}, h_{2}\not\equiv0$.
Thus, our results are new even for the second order Laplacian system;
\item[(\uppercase\expandafter{\romannumeral3})]
Because of the coupling relationship of $u$ and $v$, our proofs in the present paper become more difficult and complex than those in \cite{Qi2023}.
Especially, such difficult and complexity can be embodied in the proofs of the compactness of Palais-Smale sequence.
\end{itemize}
\par
Next, we state our main results.
\vskip2mm
 \noindent
{\bf Theorem 1.1. } {\it Assume that $(\phi_{1})$-$(\phi_{4})$, $(A)$, $(B)$ and $(H)$ hold.
Then for each $\lambda+\mu\in (0,\Lambda_{0})$ and $h_{1},h_{2}\in L^{2}(\Omega)$ with $\|h_{1}\|_{2}^{2}+\|h_{2}\|_{2}^{2}\in (0, m_{\lambda,\mu}]$, where
\begin{eqnarray}
\label{aa1}& &
\Lambda_{0}
=
\left[(2C_{0})^{-1}\left(\rho_{0}-S_{2}
\left(\epsilon_{1}^{2}+\epsilon_{2}^{2}\right)\right)\right]^{(\alpha+\beta-q)/(\alpha+\beta-2)},\\
\label{aa2}& &
m_{\lambda,\mu}
=2^{-1}S_{2}^{-1}t_{\lambda,\mu}^{2}\left(2^{-1}\rho_{0}-2^{-1}S_{2}\left(\epsilon_{1}^{2}+\epsilon_{2}^{2}\right)
-C_{0}(\lambda+\mu)^{(\alpha+\beta-2)/(\alpha+\beta-q)}\right)
\left(\max\left\{(2\epsilon_{1}^{2})^{-1},(2\epsilon_{2}^{2})^{-1}\right\}\right)^{-1},\\
\label{aa3}& &
C_{0}
=
q^{-1}S_{\alpha+\beta}^{q}
\max\left\{\|a\|_{\frac{\alpha+\beta}{\alpha+\beta-q}},\|c\|_{\frac{\alpha+\beta}{\alpha+\beta-q}}\right\}
\alpha_{0}^{(q-2)/(\alpha+\beta-q)}
+(\alpha+\beta)^{-1}\|b\|_{\infty}S_{\alpha+\beta}^{\alpha+\beta}
\alpha_{0}^{(\alpha+\beta-2)/(\alpha+\beta-q)},\nonumber\\
\label{aa4}& &
\alpha_{0}
=(2-q)(\alpha+\beta)q^{-1}(\alpha+\beta-2)^{-1}\|b\|_{\infty}^{-1}S_{\alpha+\beta}^{q-(\alpha+\beta)}
\max\left\{\|a\|_{\frac{\alpha+\beta}{\alpha+\beta-q}},\|c\|_{\frac{\alpha+\beta}{\alpha+\beta-q}}\right\},\nonumber\\
\label{aa5}& &
t_{\lambda,\mu}
=
(\lambda+\mu)^{\frac{1}{\alpha+\beta-q}}\alpha_{0}^{\frac{1}{\alpha+\beta-q}},
\;\;
\epsilon_{1},\epsilon_{2}>0,\;\;\epsilon_{1}^{2}+\epsilon_{2}^{2}<\rho_{0}/S_{2},\nonumber
\end{eqnarray}
$S_{2}$ denotes the best Sobolev constant for the imbedding $H_{0}^{1}(\Omega)\hookrightarrow L^{2}(\Omega)$,
problem (\ref{eq1}) has at least one nontrivial solution of positive energy and one nontrivial solution of negative energy.}
 \vskip2mm
Moreover, $(\phi_{4})$ can be deleted if we strengthen $(\phi_{1})$ to the following condition:
\begin{itemize}
 \item[$(\phi_1)'$] there exist two constants $0< \rho_0\leq \rho_1$ such that
 $ 0<\frac{2}{\alpha+\beta}\rho_{1}<\rho_0\leq \phi_{i}(s)\leq \rho_1$ for all $s\in [0,+\infty)$.
 \end{itemize}
 We obtain the following result.
 \vskip2mm
 \noindent
{\bf Theorem 1.2. } {\it Assume that $(\phi_1)'$, $(\phi_{2})$-$(\phi_{3})$, $(A)$, $(B)$ and $(H)$ hold.
Then for each $\lambda+\mu\in (0,\Lambda_{0})$ and $h_{1},h_{2}\in L^{2}(\Omega)$ with $\|h_{1}\|_{2}^{2}+\|h_{2}\|_{2}^{2}\in (0, m_{\lambda,\mu}]$, where $\Lambda_{0}$ and $m_{\lambda,\mu}$ are given by
(\ref{aa1}) and (\ref{aa2}), respectively,
problem (\ref{eq1}) admits at least one nontrivial solution of positive energy and one nontrivial solution of negative energy.}
\vskip2mm
 \noindent
{\bf Remark 1.1.} By $(\phi_1)$, $(A)$, $(B)$ and $(H)$, one can see that the functional $J_{\lambda,\mu}$ has the mountain pass geometry,
where $J_{\lambda,\mu}$ is the energy functional of problem (\ref{eq1}) and its definition will be given in Section 2.
Furthermore, using above conditions one more time, there exists a Palais-Smale sequence of $J_{\lambda,\mu}$ in a neighbourhood of the origin.
Then the existence of two non-trivial critical points of $J_{\lambda,\mu}$ can be shown by establishing the Palais-Smale condition.
Actually, once we could have the boundedness of Palais-Smale sequences in hand, one can expect the strong convergence of Palais-Smale sequences by $(\phi_2)$-$(\phi_3)$ and \cite[Lemma 3.5, Lemma 5.4]{Jeanjean2022}.
Hence, the existence of two non-trivial critical points of $J_{\lambda,\mu}$ is transformed into proving that any Palais-Smale sequence $\{(u_{n},v_{n})\}$ of $J_{\lambda,\mu}$ is bounded in $W$.
However, it is actually not an easy work for our problem.
The main reason is that the operator is quasi-linear, which causes a coupling term of $\Phi$ and $\phi$
in the process of obtaining the boundedness.
Although the assumption $(\phi_1)$ can be used to convert this coupling term to $\left(\frac{\rho_0}{2}-\frac{\rho_1}{\alpha+\beta}\right)\|(u_{n},v_{n})\|^{2}$, it cannot guarantee that its coefficient is a positive constant.
A straightforward way to deal with this is to make $\frac{2}{\alpha+\beta}\rho_{1}<\rho_0$ (see $(\phi_1)'$).
Moreover, inspired by the relationship $(1.4)$ of \cite{Pomponio2021}, we can also make a additional assumption  $(\phi_4)$ on $\phi$ to overcome this difficulty.
\vskip2mm
\par
The remainder of this paper focuses on some notations and preliminaries, proofs of the main results, and some examples that illustrate our results.

\section{Preliminaries}
\vskip2mm
 \par
In this section, we present some basic facts regarding Sobolev spaces and results that will be used in this work.
\par
For $1\leq p<\infty$, $ L^{p}(\Omega)$ is the usual Lebesgue space with norm $\|u\|_{p}^{p}:=\int_{\Omega}|u|^{p}dx$.
The Sobolev space $H_{0}^{1}(\Omega)$ with norm $\|u\|^{2}=\|\nabla u\|_{2}^{2}=\int_{\Omega}|\nabla u|^{2}dx$.
We denote the dual space of $H^{1}_{0}(\Omega)$ by $(H^{1}_{0}(\Omega))^{'}$.
We denote the weak convergence by $\rightharpoonup$ and denote the strong convergence by $\rightarrow$.
Finally, $S_{i}$ denotes the best Sobolev constants for the imbeddings $H_{0}^{1}(\Omega)\hookrightarrow L^{i}(\Omega)$, where $1<i<2^{\ast}=\frac{2N}{N-2}$.
\vskip2mm
 \noindent
{\bf Proposition 2.1. }(Sobolev imbedding theorem \cite{Adams2003})
Let $\Omega$ be an open subset of $\mathbb{R}^{N}$.
If $|\Omega|<\infty$, the following imbeddings are continuous:
\begin{eqnarray*}
H_{0}^{1}(\Omega)\hookrightarrow L^{\gamma}(\Omega),\;\;1\leq\gamma\leq2^{\ast}.
\end{eqnarray*}
Furthermore, the following imbeddings are compact:
\begin{eqnarray*}
H_{0}^{1}(\Omega)\hookrightarrow L^{\gamma}(\Omega),\;\;1\leq\gamma<2^{\ast},
\end{eqnarray*}
where
\begin{eqnarray*}
2^{\ast}= \begin{cases}
               \frac{2N}{N-2}, & \mbox{  if } N>2,\\
               +\infty, & \mbox{ if } N\leq2.
\end{cases}
\end{eqnarray*}
\par
Let $W=H_{0}^{1}(\Omega)\times H_{0}^{1}(\Omega)$ be the Hilbert space with norm
\begin{eqnarray*}
\|z\|
=
\left(\int_{\Omega}(|\nabla u|^{2}+|\nabla v|^{2})dx\right)^{\frac{1}{2}},
\end{eqnarray*}
where $z=(u,v)\in W$. Moreover, a pair of functions $(u,v)\in W$ is said to be a weak solution of problem (\ref{eq1}) if
\begin{eqnarray*}
& &\int_{\Omega}\left[\phi_{1}\left(\frac{u^{2}+|\nabla u|^{2}}{2}\right)(u\varphi_{1}+\nabla u\cdot \nabla \varphi_{1})
+
\phi_{2}\left(\frac{v^{2}+|\nabla v|^{2}}{2}\right)(v\varphi_{2}+\nabla v\cdot \nabla \varphi_{2} )\right]dx\\
&&
-\int_{\Omega}\left(\lambda a(x)|u|^{q-2}u\varphi_{1}+\mu c(x)|v|^{q-2}v\varphi_{2}\right)dx
-\frac{1}{\alpha+\beta}\int_{\Omega}b(x)(\alpha|u|^{\alpha-2}u|v|^{\beta}\varphi_{1}+\beta|v|^{\beta-2}v|u|^{\alpha}\varphi_{2})dx\\
&&-\int_{\Omega}h_{1}(x)\varphi_{1}dx-\int_{\Omega}h_{2}(x)\varphi_{2}dx
=0,
\;\;\forall (\varphi_{1},\varphi_{2})\in W.
\end{eqnarray*}
Furthermore, if $(u,v)\neq(0,0)$, then we call $(u,v)$ is a nontrivial solution.
The corresponding energy functional of problem (\ref{eq1}) is defined by
\begin{eqnarray}\label{2.1.1}
J_{\lambda,\mu}(u,v)
&= &
\int_{\Omega}\left[\Phi_{1}\left(\frac{u^{2}+|\nabla u|^{2}}{2}\right)+\Phi_{2}\left(\frac{v^{2}+|\nabla v|^{2}}{2}\right)\right]dx
-\frac{1}{q}\int_{\Omega}\left(\lambda a(x)|u|^{q}+\mu c(x)|v|^{q}\right)dx\nonumber\\
& &
-\frac{1}{\alpha+\beta}\int_{\Omega}b(x)|u|^{\alpha}|v|^{\beta}dx
-\int_{\Omega}h_{1}(x)udx-\int_{\Omega}h_{2}(x)vdx
\end{eqnarray}
for $(u,v)\in W$.
Assumptions $(\phi_1)$, $(H)$ and $1<q<2<\alpha+\beta<2^{\ast}$ ensure that the functional $J_{\lambda,\mu}$ is well defined in $W$ for any $\lambda>0$.
Furthermore, if $\phi_{i}\in C([0,+\infty),\mathbb{R})$ ($i=1,2$), $(\phi_1)$, $(H)$ and $1<q<2<\alpha+\beta<2^{\ast}$ hold, the functional $J_{\lambda,\mu}\in C^{1}(W,\mathbb R)$ for any $\lambda>0$, and the derivation of functional $J_{\lambda,\mu}$ satisfies
\begin{eqnarray}\label{2.1.2}
\langle J'_{\lambda,\mu}(u,v),(\varphi_{1},\varphi_{2})\rangle
&= &\int_{\Omega}\left[\phi_{1}\left(\frac{u^{2}+|\nabla u|^{2}}{2}\right)(u\varphi_{1}+\nabla u\cdot \nabla \varphi_{1})
+
\phi_{2}\left(\frac{v^{2}+|\nabla v|^{2}}{2}\right)(v\varphi_{2}+\nabla v\cdot \nabla \varphi_{2} )\right]dx\nonumber\\
& &
-\int_{\Omega}\left(\lambda a(x)|u|^{q-2}u\varphi_{1}+\mu c(x)|v|^{q-2}v\varphi_{2}\right)dx
-\int_{\Omega}h_{1}(x)\varphi_{1}dx-\int_{\Omega}h_{2}(x)\varphi_{2}dx\nonumber\\
& &
-\frac{1}{\alpha+\beta}\int_{\Omega}b(x)(\alpha|u|^{\alpha-2}u|v|^{\beta}\varphi_{1}
+\beta|v|^{\beta-2}v|u|^{\alpha}\varphi_{2})dx
\end{eqnarray}
for any $(u,v), (\varphi_{1},\varphi_{2})\in W$.
It is easy to see that (\ref{eq1}) is the Euler-Lagrange equation of the functional $J_{\lambda,\mu}$.
Hence, finding a weak solution for the problem (\ref{eq1}) is equivalent to finding a critical point for the functional $J_{\lambda,\mu}$.
\par
Assume that $\varphi \in C^{1}(X,\mathbb{R})$. An sequence $\{u_{n}\}$ is called as the Palais-Smale sequence of $\varphi$ if $\varphi(u_{n})$ is bounded for all $n\in N$ and $\varphi'(u_{n})\rightarrow0$ as $n\rightarrow\infty$. If any Palais-Smale sequence $\{u_{n}\}$ of $\varphi$ has a convergent subsequence, we call that $\varphi$ satisfies the
Palais-Smale condition ((PS)-condition for short).
\vskip2mm
 \noindent
 {\bf Lemma 2.1.} ( Ekeland's variational principle \cite{Mawhin1989})
 Let $X$ be a complete metric space with metric $d$ and $\varphi:X\rightarrow\mathbb{R}$ be a lower semicontinuous function, bounded from below and not identical to $+\infty$. Let $\varepsilon>0$ be given and $U\in X$ such that
 \begin{eqnarray*}
\varphi(U)\leq \inf_{M}\varphi+\varepsilon.
\end{eqnarray*}
Then there exists $V\in X$ such that
\begin{eqnarray*}
\varphi(V)\leq \varphi(U),\;\;d(U,V)\leq1,
\end{eqnarray*}
and for each $E\in X$, one has
\begin{eqnarray*}
\varphi(V)\leq \varphi(E)+\varepsilon d(V,E).
\end{eqnarray*}
\par
By the Ekeland's variational principle, it is easy to obtain the following corollary.
\vskip2mm
 \noindent
 {\bf Lemma 2.2.} (\cite{Mawhin1989})
 Suppose that $X$ is a Banach space, $M\subset X$ is closed, $\varphi\in C^{1}(X,\mathbb{R})$ is bounded from below on $M$ and satisfies the (PS)-condition. Then $\varphi$ attains its infimum on $M$.
 \vskip2mm
 \noindent
 {\bf Lemma 2.3.} (Mountain pass theorem \cite{Rabinowitz1986}) Let $X$ be a real Banach space and $\varphi \in C^{1}(X,\mathbb{R})$, $\varphi(0)=0$ satisfy (PS)-condition. Suppose that $\varphi$ satisfies the following conditions:
\begin{itemize}
\item[$(i)$] there exists a constant $\rho>0$ and $\alpha>0$ such that $\varphi|_{\partial B_{\rho}(0)}\geq\alpha$, where $B_{\rho}=\{\omega\in X:\|\omega\|_{X}<\rho\}$;
\item[$(ii)$] there exists $\omega\in X\setminus\bar{B}_{\rho}(0)$ such that $\varphi(\omega)\leq0$.
\end{itemize}
 Then $\varphi$ has a critical value $c_{\ast}\geq\alpha$ with
 \begin{eqnarray*}
 c_{\ast}:=\inf_{\gamma\in\Gamma}\max_{t\in[0,1]}\varphi(\gamma(t)),
\end{eqnarray*}
where
\begin{eqnarray*}
\Gamma:=\{\gamma\in C([0,1],X):\gamma(0)=0, \; \gamma(1)= \omega\}.
\end{eqnarray*}

\section{Proofs of the main results}
\vskip2mm
 \noindent
 {\bf Lemma 3.1.} Suppose that $(\phi_1)$, $(A)$ and $(B)$ hold. Then the following conclusions hold:
 \begin{itemize}
\item[$(i)$]
$\displaystyle{
\frac{\rho_{0}}{2}\|(u,v)\|^{2}
\leq
\int_{\Omega}\left[\Phi_{1}\left(\frac{u^{2}+|\nabla u|^{2}}{2}\right)
+\Phi_{2}\left(\frac{v^{2}+|\nabla v|^{2}}{2}\right)\right]dx
\leq
\frac{\rho_{1}}{2}\|(u,v)\|^{2}}$,
\item[$(ii)$]
$\displaystyle{\int_{\Omega}(\lambda a(x)|u|^{q}+\mu c(x)|v|^{q})dx
\leq (\lambda+\mu)S_{\alpha+\beta}^{q}
\max\left\{\|a\|_{\frac{\alpha+\beta}{\alpha+\beta-q}},\|c\|_{\frac{\alpha+\beta}{\alpha+\beta-q}}\right\}
\|(u,v)\|^{q}}$,
\item[$(iii)$]
$\displaystyle{\int_{\Omega}b(x)|u|^{\alpha}|v|^{\beta}dx
\leq
\|b\|_{\infty}S_{\alpha+\beta}^{\alpha+\beta}\|(u,v)\|^{\alpha+\beta}}$
\end{itemize}
for any $(u,v)\in W$.
\vskip2mm
\noindent
{\bf Proof.} Initially, we shall prove the item $(i)$.
By $(\phi_1)$, we deduce that
\begin{eqnarray}\label{3.1.1}
\int_{\Omega}\left[\Phi_{1}\left(\frac{u^{2}+|\nabla u|^{2}}{2}\right)+\Phi_{2}\left(\frac{v^{2}+|\nabla v|^{2}}{2}\right)\right]dx
& \geq &
 \int_{\Omega}\left[\rho_{0}\frac{u^{2}+|\nabla u|^{2}}{2}+\rho_{0}\frac{v^{2}+|\nabla v|^{2}}{2}\right]dx\nonumber\\
& \geq &
 \frac{\rho_{0}}{2}\int_{\Omega}(|\nabla u|^{2}+|\nabla v|^{2})dx\nonumber\\
& = &
\frac{\rho_{0}}{2}\|(u,v)\|^{2}
\end{eqnarray}
and
\begin{eqnarray}\label{3.1.2}
\int_{\Omega}\left[\Phi_{1}\left(\frac{u^{2}+|\nabla u|^{2}}{2}\right)+\Phi_{2}\left(\frac{v^{2}+|\nabla v|^{2}}{2}\right)\right]dx
& \leq &
 \int_{\Omega}\left[\rho_{1}\frac{u^{2}+|\nabla u|^{2}}{2}+\rho_{1}\frac{v^{2}+|\nabla v|^{2}}{2}\right]dx\nonumber\\
& \leq &
 \frac{\rho_{1}}{2}\int_{\Omega}(|\nabla u|^{2}+|\nabla v|^{2})dx\nonumber\\
& = &
\frac{\rho_{1}}{2}\|(u,v)\|^{2}.
\end{eqnarray}
This ends the proof of the item $(i)$.
Next we shall prove the item $(ii)$.
Note that $a\in L^{\frac{\alpha+\beta}{\alpha+\beta-q}+\tau_{1}}(\Omega)$,
$c \in L^{\frac{\alpha+\beta}{\alpha+\beta-q}+\tau_{2}}(\Omega)$ implies that
$a\in L^{\frac{\alpha+\beta}{\alpha+\beta-q}}(\Omega)$,
 $c \in L^{\frac{\alpha+\beta}{\alpha+\beta-q}}(\Omega)$.
Then, by H\"older inequality and Proposition 2.1, we have
\begin{eqnarray}\label{3.1.3}
&  &
\lambda\int_{\Omega}a(x)|u|^{q}dx+\mu\int_{\Omega} c(x)|v|^{q}dx\nonumber\\
& \leq &
\lambda\int_{\Omega}|a(x)||u|^{q}dx
+
\mu\int_{\Omega}| c(x)||v|^{q}dx\nonumber\\
& \leq &
\lambda\left(\int_{\Omega}|a(x)|^{\frac{\alpha+\beta}{\alpha+\beta-q}}dx\right)^{\frac{\alpha+\beta-q}{\alpha+\beta}}
\left(\int_{\Omega}|u|^{q\cdot\frac{\alpha+\beta}{q}}dx\right)^{\frac{q}{\alpha+\beta}}
+
\mu\left(\int_{\Omega}|c(x)|^{\frac{\alpha+\beta}{\alpha+\beta-q}}dx\right)^{\frac{\alpha+\beta-q}{\alpha+\beta}}
\left(\int_{\Omega}|v|^{q\cdot\frac{\alpha+\beta}{q}}dx\right)^{\frac{q}{\alpha+\beta}}\nonumber\\
& =  &
\lambda\|a\|_{\frac{\alpha+\beta}{\alpha+\beta-q}}\|u\|_{\alpha+\beta}^{q}
+
\mu\|c\|_{\frac{\alpha+\beta}{\alpha+\beta-q}}\|v\|_{\alpha+\beta}^{q}\nonumber\\
& \leq &
S_{\alpha+\beta}^{q}\left(\lambda\|a\|_{\frac{\alpha+\beta}{\alpha+\beta-q}}\|\nabla u\|_{2}^{q}
+
\mu\|c\|_{\frac{\alpha+\beta}{\alpha+\beta-q}}\|\nabla v\|_{2}^{q}\right)\nonumber\\
& \leq &
S_{\alpha+\beta}^{q}\left(\lambda\|a\|_{\frac{\alpha+\beta}{\alpha+\beta-q}}
+\mu\|c\|_{\frac{\alpha+\beta}{\alpha+\beta-q}}\right)
\|(u,v)\|^{q}\nonumber\\
& \leq &
(\lambda+\mu)S_{\alpha+\beta}^{q}
\max\left\{\|a\|_{\frac{\alpha+\beta}{\alpha+\beta-q}},\|c\|_{\frac{\alpha+\beta}{\alpha+\beta-q}}\right\}\|(u,v)\|^{q}.
\end{eqnarray}
Finally, we shall prove the item $(iii)$.
By $(B)$, the Young's inequality and Proposition 2.1, we mention that
\begin{eqnarray}\label{3.1.4}
\int_{\Omega}b(x)|u|^{\alpha}|v|^{\beta}dx
& \leq &
\int_{\Omega}|b(x)||u|^{\alpha}|v|^{\beta}dx\nonumber\\
& \leq &
\|b\|_{\infty}\int_{\Omega}|u|^{\alpha}|v|^{\beta}dx\nonumber\\
& \leq &
\|b\|_{\infty}\int_{\Omega}\left(\frac{\alpha}{\alpha+\beta}|u|^{\alpha\cdot\frac{\alpha+\beta}{\alpha}}
+\frac{\beta}{\alpha+\beta}|v|^{\beta\cdot\frac{\alpha+\beta}{\beta}}\right)dx\nonumber\\
&= &
\|b\|_{\infty}\left(\frac{\alpha}{\alpha+\beta}\|u\|^{\alpha+\beta}_{\alpha+\beta}
+\frac{\beta}{\alpha+\beta}\|v\|^{\alpha+\beta}_{\alpha+\beta}\right)\nonumber\\
& \leq &
\|b\|_{\infty}S_{\alpha+\beta}^{\alpha+\beta}\left(\frac{\alpha}{\alpha+\beta}\|\nabla u\|^{\alpha+\beta}_{2}
+\frac{\beta}{\alpha+\beta}\|v\|^{\alpha+\beta}_{2}\right)\nonumber\\
& \leq &
\|b\|_{\infty}S_{\alpha+\beta}^{\alpha+\beta}\|(u,v)\|^{\alpha+\beta}.
\end{eqnarray}
The proof for this lemma now is completed.
\qed
 \vskip2mm
 \noindent
 {\bf Lemma 3.2.} Assume that $(\phi_1)$, $(A)$, $(B)$ and $(H)$ hold. Then for each $\lambda+\mu\in (0,\Lambda_{0})$ and $h_{1},h_{2}\in L^{2}(\Omega)$ with $\|h_{1}\|_{2}^{2}+\|h_{2}\|_{2}^{2}\in (0, m_{\lambda,\mu}]$, where
$\Lambda_{0}$ and $m_{\lambda,\mu}$ were given in (\ref{aa1}) and (\ref{aa2}) respectively,
there exists a positive constant $\rho_{\lambda,\mu}$ such that $J_{\lambda,\mu}(u,v)>0$ whenever $\|(u,v)\|=\rho_{\lambda,\mu}$.
 \vskip2mm
  \noindent
 {\bf Proof.} Let $\epsilon_{1},\epsilon_{2}>0$ with $\epsilon_{1}^{2}+\epsilon_{2}^{2}<\rho_{0}/S_{2}$.
 It follows from $(H)$, H\"older inequality, Proposition 2.1 and the Young's inequality that
\begin{eqnarray}\label{3.2.1}
\int_{\Omega}|h_{1}||u|dx
\leq S_{2}\|h_{1}\|_{2}\|\nabla u\|_{2}
\leq S_{2}\left(\frac{\epsilon_{1}^{2}}{2}\|\nabla u\|_{2}^{2}+\frac{1}{2\epsilon_{1}^{2}}\|h_{1}\|_{2}^{2}\right)
\leq S_{2}\left(\frac{\epsilon_{1}^{2}}{2}\|(u,v)\|^{2}+\frac{1}{2\epsilon_{1}^{2}}\|h_{1}\|_{2}^{2}\right).
\end{eqnarray}
Similarly, we also have
\begin{eqnarray}\label{3.2.2}
\int_{\Omega}|h_{2}||v|dx
\leq S_{2}\left(\frac{\epsilon_{2}^{2}}{2}\|(u,v)\|^{2}+\frac{1}{2\epsilon_{2}^{2}}\|h_{2}\|_{2}^{2}\right).
\end{eqnarray}
Then, by (\ref{2.1.1}), (\ref{3.1.1}), (\ref{3.1.3}), (\ref{3.1.4}), (\ref{3.2.1}) and (\ref{3.2.2}), we mention that
\begin{eqnarray}\label{3.2.3}
& &J_{\lambda,\mu}(u,v)\nonumber\\
& =&
\int_{\Omega}\left[\Phi_{1}\left(\frac{u^{2}+|\nabla u|^{2}}{2}\right)+\Phi_{2}\left(\frac{v^{2}+|\nabla v|^{2}}{2}\right)\right]dx
-
\frac{1}{q}\int_{\Omega}\left(\lambda a(x)|u|^{q}+\mu c(x)|v|^{q}\right)dx\nonumber\\
& &
-\frac{1}{\alpha+\beta}\int_{\Omega}b(x)|u|^{\alpha}|v|^{\beta}dx
-\int_{\Omega}h_{1}udx-\int_{\Omega}h_{2}vdx\nonumber\\
& \geq&
\frac{\rho_{0}}{2}\|(u,v)\|^{2}
-
(\lambda+\mu)\frac{1}{q}S_{\alpha+\beta}^{q}
\max\left\{\|a\|_{\frac{\alpha+\beta}{\alpha+\beta-q}},\|c\|_{\frac{\alpha+\beta}{\alpha+\beta-q}}\right\}
\|(u,v)\|^{q}
-
\frac{1}{\alpha+\beta}\|b\|_{\infty}S_{\alpha+\beta}^{\alpha+\beta}\|(u,v)\|^{\alpha+\beta}\nonumber\\
& &
-\frac{S_{2}}{2}\left(\epsilon_{1}^{2}+\epsilon_{2}^{2}\right)\|(u,v)\|^{2}
-\left(\frac{S_{2}}{2\epsilon_{1}^{2}}\|h_{1}\|_{2}^{2}
+\frac{S_{2}}{2\epsilon_{2}^{2}}\|h_{2}\|_{2}^{2}\right)\nonumber\\
& \geq&
\left(\frac{\rho_{0}}{2}-\frac{S_{2}}{2}\left(\epsilon_{1}^{2}+\epsilon_{2}^{2}\right)\right)\|(u,v)\|^{2}
-
(\lambda+\mu)\frac{1}{q}S_{\alpha+\beta}^{q}
\max\left\{\|a\|_{\frac{\alpha+\beta}{\alpha+\beta-q}},\|c\|_{\frac{\alpha+\beta}{\alpha+\beta-q}}\right\}
\|(u,v)\|^{q}\nonumber\\
& &
-\frac{1}{\alpha+\beta}\|b\|_{\infty}S_{\alpha+\beta}^{\alpha+\beta}\|(u,v)\|^{\alpha+\beta}
-\max\left\{\frac{S_{2}}{2\epsilon_{1}^{2}},\frac{S_{2}}{2\epsilon_{2}^{2}}\right\}
\left(\|h_{1}\|_{2}^{2}+\|h_{2}\|_{2}^{2}\right)\nonumber\\
& =&
\|(u,v)\|^{2}
\Bigg\{\frac{\rho_{0}}{2}-\frac{S_{2}}{2}\left(\epsilon_{1}^{2}+\epsilon_{2}^{2}\right)
-(\lambda+\mu)\frac{1}{q}S_{\alpha+\beta}^{q}
\max\left\{\|a\|_{\frac{\alpha+\beta}{\alpha+\beta-q}},\|c\|_{\frac{\alpha+\beta}{\alpha+\beta-q}}\right\}
\|(u,v)\|^{q-2}\nonumber\\
& &
-\frac{1}{\alpha+\beta}\|b\|_{\infty}S_{\alpha+\beta}^{\alpha+\beta}\|(u,v)\|^{\alpha+\beta-2}\Bigg\}
-\max\left\{\frac{S_{2}}{2\epsilon_{1}^{2}},\frac{S_{2}}{2\epsilon_{2}^{2}}\right\}
\left(\|h_{1}\|_{2}^{2}+\|h_{2}\|_{2}^{2}\right)
\end{eqnarray}
for $(u,v)\in W$.
Let
\begin{eqnarray*}
f(t)
=
(\lambda+\mu)\frac{1}{q}S_{\alpha+\beta}^{q}
\max\left\{\|a\|_{\frac{\alpha+\beta}{\alpha+\beta-q}},\|c\|_{\frac{\alpha+\beta}{\alpha+\beta-q}}\right\}
t^{q-2}
+
\frac{1}{\alpha+\beta}\|b\|_{\infty}S_{\alpha+\beta}^{\alpha+\beta}t^{\alpha+\beta-2},
\;\;
t\in(0,+\infty).
\end{eqnarray*}
It follows from $f'(t_{\lambda,\mu})=0$ that
\begin{eqnarray*}
t_{\lambda,\mu}
=
(\lambda+\mu)^{\frac{1}{\alpha+\beta-q}}
\left(\frac{(2-q)(\alpha+\beta)
\max\left\{\|a\|_{\frac{\alpha+\beta}{\alpha+\beta-q}},\|c\|_{\frac{\alpha+\beta}{\alpha+\beta-q}}\right\}}
{q(\alpha+\beta-2)\|b\|_{\infty}S_{\alpha+\beta}^{\alpha+\beta-q}}\right)^{\frac{1}{\alpha+\beta-q}}>0,
\end{eqnarray*}
which is a unique stationary point in $(0,+\infty)$.
By $1<q<2<\alpha+\beta$, we also mention that
\begin{eqnarray*}
\lim_{t\rightarrow0^{+}}f(t)=\lim_{t\rightarrow+\infty }f(t)=+\infty.
\end{eqnarray*}
Thus, $t_{\lambda,\mu}$ is a minimum point of $f(t)$ and
\begin{eqnarray}\label{3.2.5}
f(t_{\lambda,\mu})
=
(\lambda+\mu)\frac{1}{q}S_{\alpha+\beta}^{q}
\max\left\{\|a\|_{\frac{\alpha+\beta}{\alpha+\beta-q}},\|c\|_{\frac{\alpha+\beta}{\alpha+\beta-q}}\right\}
t_{\lambda,\mu}^{q-2}
+
\frac{1}{\alpha+\beta}\|b\|_{\infty}S_{\alpha+\beta}^{\alpha+\beta}t_{\lambda,\mu}^{\alpha+\beta-2}
=
C_{0}(\lambda+\mu)^{\frac{\alpha+\beta-2}{\alpha+\beta-q}},
\end{eqnarray}
where
\begin{eqnarray*}
& &
C_{0}
=
\frac{1}{q}S_{\alpha+\beta}^{q}
\max\left\{\|a\|_{\frac{\alpha+\beta}{\alpha+\beta-q}},\|c\|_{\frac{\alpha+\beta}{\alpha+\beta-q}}\right\}
\alpha_{0}^{\frac{q-2}{\alpha+\beta-q}}
+\frac{1}{\alpha+\beta}\|b\|_{\infty}S_{\alpha+\beta}^{\alpha+\beta}
\alpha_{0}^{\frac{\alpha+\beta-2}{\alpha+\beta-q}},\\
& &
\alpha_{0}
=\frac{(2-q)(\alpha+\beta)
\max\left\{\|a\|_{\frac{\alpha+\beta}{\alpha+\beta-q}},\|c\|_{\frac{\alpha+\beta}{\alpha+\beta-q}}\right\}}
{q(\alpha+\beta-2)\|b\|_{\infty}S_{\alpha+\beta}^{\alpha+\beta-q}}.
\end{eqnarray*}
Let
\begin{eqnarray*}
\Lambda_{0}
=
\left\{\frac{1}{C_{0}}\left(\frac{\rho_{0}}{2}-\frac{S_{2}}{2}
\left(\epsilon_{1}^{2}+\epsilon_{2}^{2}\right)\right)\right\}^{\frac{\alpha+\beta-q}{\alpha+\beta-2}},
\;\;
\alpha_{\lambda,\mu}
=\max\left\{\frac{S_{2}}{2\epsilon_{1}^{2}},\frac{S_{2}}{2\epsilon_{2}^{2}}\right\}m_{\lambda,\mu}
\end{eqnarray*}
and
\begin{eqnarray*}
m_{\lambda,\mu}
=\frac{t_{\lambda,\mu}^{2}\left(\frac{\rho_{0}}{2}-\frac{S_{2}}{2}\left(\epsilon_{1}^{2}+\epsilon_{2}^{2}\right)
-C_{0}(\lambda+\mu)^{\frac{\alpha+\beta-2}{\alpha+\beta-q}}\right)}
{2\max\left\{\frac{S_{2}}{2\epsilon_{1}^{2}},\frac{S_{2}}{2\epsilon_{2}^{2}}\right\}}.
\end{eqnarray*}
Clearly, $\Lambda_{0}>0$ since $\epsilon_{1}^{2}+\epsilon_{2}^{2}<\rho_{0}/S_{2}$.
If $\lambda+\mu\in (0,\Lambda_{0})$, then $m_{\lambda,\mu}>0$ and $\alpha_{\lambda,\mu}>0$.
Hence, taking $\rho_{\lambda,\mu}=t_{\lambda,\mu}$.
For each $\lambda+\mu\in (0,\Lambda_{0})$ and $h_{1},h_{2}\in L^{2}(\Omega)$ with
$\|h_{1}\|_{2}^{2}+\|h_{2}\|_{2}^{2}\in (0, m_{\lambda,\mu}]$,
by (\ref{3.2.3}) and (\ref{3.2.5}), we have
\begin{eqnarray*}
& &J_{\lambda,\mu}(u,v)\nonumber\\
& \geq&
t_{\lambda,\mu}^{2}
\Bigg\{\frac{\rho_{0}}{2}-\frac{S_{2}}{2}\left(\epsilon_{1}^{2}+\epsilon_{2}^{2}\right)
-
(\lambda+\mu)\frac{1}{q}S_{\alpha+\beta}^{q}
\max\left\{\|a\|_{\frac{\alpha+\beta}{\alpha+\beta-q}},\|c\|_{\frac{\alpha+\beta}{\alpha+\beta-q}}\right\}
t_{\lambda,\mu}^{q-2}
-\frac{1}{\alpha+\beta}\|b\|_{\infty}S_{\alpha+\beta}^{\alpha+\beta}t_{\lambda,\mu}^{\alpha+\beta-2}\Bigg\}\nonumber\\
& &
-\max\left\{\frac{S_{2}}{2\epsilon_{1}^{2}},\frac{S_{2}}{2\epsilon_{2}^{2}}\right\}
\left(\|h_{1}\|_{2}^{2}+\|h_{2}\|_{2}^{2}\right)\nonumber\\
& =&
t_{\lambda,\mu}^{2}
\Bigg(\frac{\rho_{0}}{2}-\frac{S_{2}}{2}\left(\epsilon_{1}^{2}+\epsilon_{2}^{2}\right)
-C_{0}(\lambda+\mu)^{\frac{\alpha+\beta-2}{\alpha+\beta-q}}\Bigg)
-\max\left\{\frac{S_{2}}{2\epsilon_{1}^{2}},\frac{S_{2}}{2\epsilon_{2}^{2}}\right\}
\left(\|h_{1}\|_{2}^{2}+\|h_{2}\|_{2}^{2}\right)\nonumber\\
& =&
2\max\left\{\frac{S_{2}}{2\epsilon_{1}^{2}},\frac{S_{2}}{2\epsilon_{2}^{2}}\right\}m_{\lambda,\mu}
-\max\left\{\frac{S_{2}}{2\epsilon_{1}^{2}},\frac{S_{2}}{2\epsilon_{2}^{2}}\right\}
\left(\|h_{1}\|_{2}^{2}+\|h_{2}\|_{2}^{2}\right)\nonumber\\
& \geq&
2\max\left\{\frac{S_{2}}{2\epsilon_{1}^{2}},\frac{S_{2}}{2\epsilon_{2}^{2}}\right\}m_{\lambda,\mu}
-\max\left\{\frac{S_{2}}{2\epsilon_{1}^{2}},\frac{S_{2}}{2\epsilon_{2}^{2}}\right\}m_{\lambda,\mu}\nonumber\\
& =&
\max\left\{\frac{S_{2}}{2\epsilon_{1}^{2}},\frac{S_{2}}{2\epsilon_{2}^{2}}\right\}m_{\lambda,\mu}\nonumber\\
& =&
\alpha_{\lambda,\mu}>0
\end{eqnarray*}
for any $(u,v)\in W$ with $\|(u,v)\|=\rho_{\lambda,\mu}$.
The proof for this lemma now complete.
\qed
\vskip2mm
 \noindent
{\bf Remark 3.1.}
If $b(x)$ are sign-changing on the whole bounded domain $\Omega$,
we can also obtain that the existence of negative energy solution for system (\ref{eq1}) by Ekeland's variational principle.
However, to obtain the fact that the energy functional has a mountain pass geometry may be challenging.
The main reason is that the sign of the integral $\int_{\Omega}b(x)|\varphi_{1}|^{\alpha}|\varphi_{2}|^{\beta}dx$ is indefinite for any $(\varphi_{1},\varphi_{2})\in C_{0}^{\infty}(\Omega)\times C_{0}^{\infty}(\Omega)$.
A straightforward way to deal with this is to make $b(x)$ nonnegative on the entire bounded domain $\Omega$
(for example \cite{Huang2013}).
In the present work,
with the help of the truncation technique in \cite{Chen2013}, we just require that $b(x)$ is nonnegative only on a non-empty open subset of $\Omega$.
\vskip2mm
 \noindent
{\bf Lemma 3.3.} Assume that $(\phi_1)$ and $(B)$ hold.
Then there exist a function $(\omega_{1},\omega_{2})\in W$ with $\|(\omega_{1},\omega_{2})\|>\rho_{\lambda,\mu}$ such that $J_{\lambda,\mu}(\omega_{1},\omega_{2})<0$.
 \vskip2mm
  \noindent
{\bf Proof.} Let $\Omega_{0}\subset\Omega_{1}$ be a bounded domain, where $\Omega_{1}$ is given in $(H)$.
Choose $(\varphi_{1},\varphi_{2})\in C_{0}^{\infty}(\Omega_{0})\times C_{0}^{\infty}(\Omega_{0})$, $\varphi_{1},\varphi_{2}\geq0$, $\varphi_{1},\varphi_{2}\not\equiv0$ in $\Omega_{0}$.
Let $\varphi_{1}(x)\cdot\varphi_{2}(x)=0$, for all $x\in\Omega^{c}_{0}=\Omega\backslash \overline{\Omega}_{0}$.
Then by (\ref{3.1.2}), we have
\begin{eqnarray*}
& &
J_{\lambda,\mu}(t\varphi_{1},t\varphi_{2})\\
& \leq&
t^{2}\frac{\rho_{1}}{2}\|(\varphi_{1},\varphi_{2})\|^{2}
-
t^{q}\frac{1}{q}\int_{\Omega}(\lambda a(x)|\varphi_{1}|^{q}+\mu c(x)|\varphi_{2}|^{q})dx
-
t^{\alpha+\beta}\frac{1}{\alpha+\beta}\int_{\Omega}b(x)|\varphi_{1}|^{\alpha}|\varphi_{2}|^{\beta}dx\\
& &
-t\int_{\Omega}h_{1}\varphi_{1}dx-t\int_{\Omega}h_{2}\varphi_{2}dx\\
& =&
t^{2}\frac{\rho_{1}}{2}\|(\varphi_{1},\varphi_{2})\|^{2}
-
t^{q}\frac{1}{q}\int_{\Omega}(\lambda a(x)|\varphi_{1}|^{q}+\mu c(x)|\varphi_{2}|^{q})dx
-
t^{\alpha+\beta}\frac{1}{\alpha+\beta}\int_{\Omega_{0}}b(x)|\varphi_{1}|^{\alpha}|\varphi_{2}|^{\beta}dx\\
& &
-t\int_{\Omega}h_{1}\varphi_{1}dx-t\int_{\Omega}h_{2}\varphi_{2}dx\\
\end{eqnarray*}
and $J_{\lambda,\mu}(t\varphi_{1},t\varphi_{2})\rightarrow-\infty$ as $t\rightarrow+\infty$ since $1<q<2<\alpha+\beta$.
Therefore, there exists $t_{1}>0$ large enough such that $\|(t_{1}\varphi_{1},t_{1}\varphi_{2})\|>\rho_{\lambda,\mu}$ and $J_{\lambda,\mu}(t_{1}\varphi_{1},t_{1}\varphi_{2})<0$.
The proof for this lemma now complete.
\qed

\vskip2mm
 \noindent
{\bf Lemma 3.4.} Assume that $(\phi_1)$-$(\phi_4)$, $(A)$, $(B)$ and $(H)$ hold. Then $J_{\lambda,\mu}(u,v)$ defined by (\ref{2.1.1}) satisfies (PS) condition on $W$.
\vskip2mm
  \noindent
{\bf Proof.} Let $c\in\mathbb{R}$ and $\{(u_{n},v_{n})\}$ be a $(PS)_{c}$ sequence of $J_{\lambda,\mu}$ in $W$, that is,
\begin{eqnarray}\label{3.4.1}
\lim_{n\rightarrow\infty}J_{\lambda,\mu}(u_{n},v_{n})=c
\;\;\mbox{and}\;\;
\lim_{n\rightarrow\infty}J'_{\lambda,\mu}(u_{n},v_{n})=0.
\end{eqnarray}
We first claim that $\{(u_{n},v_{n})\}$ is bounded in $W$.
Clearly, by (\ref{3.4.1}), we have
\begin{eqnarray*}
J_{\lambda,\mu}(u_{n},v_{n})-\frac{1}{\alpha+\beta}\langle J'_{\lambda,\mu}(u_{n},v_{n}),(u_{n},v_{n}) \rangle
\rightarrow
c
<
c+1,\;\;\mbox{as}\;\; n\rightarrow\infty.
\end{eqnarray*}
Thus, there exists $n_{0}>0$ large enough such that
\begin{eqnarray}\label{3.4.2}
J_{\lambda,\mu}(u_{n},v_{n})-\frac{1}{\alpha+\beta}\langle J'_{\lambda,\mu}(u_{n},v_{n}),(u_{n},v_{n}) \rangle
<
c+1
\end{eqnarray}
for any $n>n_{0}$.
In addition, by $(\phi_1)$, $(\phi_4)$, $(A)$, $(H)$, (\ref{2.1.1}), (\ref{3.1.1}) and (\ref{3.1.3}), we have
\begin{eqnarray}\label{3.4.3}
& &
J_{\lambda,\mu}(u_{n},v_{n})-\frac{1}{\alpha+\beta}\langle J'_{\lambda,\mu}(u_{n},v_{n}),(u_{n},v_{n}) \rangle\nonumber\\
&= &
\int_{\Omega}\left[\Phi_{1}\left(\frac{u_{n}^{2}+|\nabla u_{n}|^{2}}{2}\right)
-\frac{1}{\alpha+\beta}\phi_{1}\left(\frac{u_{n}^{2}+|\nabla u_{n}|^{2}}{2}\right)(u_{n}^{2}+|\nabla u_{n}|^{2})\right]dx\nonumber\\
& &
+\int_{\Omega}\left[\Phi_{2}\left(\frac{v_{n}^{2}+|\nabla v_{n}|^{2}}{2}\right)
-\frac{1}{\alpha+\beta}\phi_{2}\left(\frac{v_{n}^{2}+|\nabla v_{n}|^{2}}{2}\right)(v_{n}^{2}+|\nabla v_{n}|^{2})\right]dx\nonumber\\
& &
-\left(\frac{1}{q}-\frac{1}{\alpha+\beta}\right)\int_{\Omega}(\lambda a(x)|u_{n}|^{q}+\mu c(x)|v_{n}|^{q})dx
-\left(1-\frac{1}{\alpha+\beta}\right)\int_{\Omega}(h_{1}u_{n}+h_{2}v_{n})dx
\nonumber\\
& \geq&
\int_{\Omega}\left[\phi_{1}\left(\frac{u_{n}^{2}+|\nabla u_{n}|^{2}}{2}\right)\frac{u_{n}^{2}+|\nabla u_{n}|^{2}}{2}  -\frac{1}{\alpha+\beta}\phi_{1}\left(\frac{u_{n}^{2}+|\nabla u_{n}|^{2}}{2}\right)(u_{n}^{2}+|\nabla u_{n}|^{2})\right]dx\nonumber\\
& &
+\int_{\Omega}\left[\phi_{2}\left(\frac{v_{n}^{2}+|\nabla v_{n}|^{2}}{2}\right)\frac{v_{n}^{2}+|\nabla v_{n}|^{2}}{2}  -\frac{1}{\alpha+\beta}\phi_{2}\left(\frac{v_{n}^{2}+|\nabla v_{n}|^{2}}{2}\right)(v_{n}^{2}+|\nabla v_{n}|^{2})\right]dx\nonumber\\
& &
-\left(\frac{1}{q}-\frac{1}{\alpha+\beta}\right)(\lambda+\mu)S_{\alpha+\beta}^{q}
\max\left\{\|a\|_{\frac{\alpha+\beta}{\alpha+\beta-q}},\|c\|_{\frac{\alpha+\beta}{\alpha+\beta-q}}\right\}
\|(u_{n},v_{n})\|^{q}\nonumber\\
& &
-\left(1-\frac{1}{\alpha+\beta}\right)S_{2}(\|h_{1}\|_{2}+\|h_{2}\|_{2})\|(u_{n},v_{n})\|
\nonumber\\
&=&
\frac{\alpha+\beta-2}{\alpha+\beta}
\int_{\Omega}\left[\phi_{1}\left(\frac{u_{n}^{2}+|\nabla u_{n}|^{2}}{2}\right)
\frac{u_{n}^{2}+|\nabla u_{n}|^{2}}{2}
+\phi_{2}\left(\frac{v_{n}^{2}+|\nabla v_{n}|^{2}}{2}\right)
\frac{v_{n}^{2}+|\nabla v_{n}|^{2}}{2}
\right]dx\nonumber\\
& &
-\left(\frac{1}{q}-\frac{1}{\alpha+\beta}\right)(\lambda+\mu)S_{\alpha+\beta}^{q}
\max\left\{\|a\|_{\frac{\alpha+\beta}{\alpha+\beta-q}},\|c\|_{\frac{\alpha+\beta}{\alpha+\beta-q}}\right\}
\|(u_{n},v_{n})\|^{q}\nonumber\\
& &
-\left(1-\frac{1}{\alpha+\beta}\right)S_{2}(\|h_{1}\|_{2}+\|h_{2}\|_{2})\|(u_{n},v_{n})\|
\nonumber\\
& \geq&
\frac{(\alpha+\beta-2)\rho_{0}}{2(\alpha+\beta)}
\int_{\Omega}(u_{n}^{2}+|\nabla u_{n}|^{2}+v_{n}^{2}+|\nabla v_{n}|^{2})dx
-\left(1-\frac{1}{\alpha+\beta}\right)S_{2}(\|h_{1}\|_{2}+\|h_{2}\|_{2})\|(u_{n},v_{n})\|\nonumber\\
& &
-\left(\frac{1}{q}-\frac{1}{\alpha+\beta}\right)(\lambda+\mu)S_{\alpha+\beta}^{q}
\max\left\{\|a\|_{\frac{\alpha+\beta}{\alpha+\beta-q}},\|c\|_{\frac{\alpha+\beta}{\alpha+\beta-q}}\right\}
\|(u_{n},v_{n})\|^{q}
\nonumber\\
& \geq&
\frac{(\alpha+\beta-2)\rho_{0}}{2(\alpha+\beta)}
\int_{\Omega}(|\nabla u_{n}|^{2}+|\nabla v_{n}|^{2})dx
-\left(\frac{1}{q}-\frac{1}{\alpha+\beta}\right)(\lambda+\mu)S_{\alpha+\beta}^{q}
\max\left\{\|a\|_{\frac{\alpha+\beta}{\alpha+\beta-q}},\|c\|_{\frac{\alpha+\beta}{\alpha+\beta-q}}\right\}
\|(u_{n},v_{n})\|^{q}\nonumber\\
& &
-\left(1-\frac{1}{\alpha+\beta}\right)S_{2}(\|h_{1}\|_{2}+\|h_{2}\|_{2})\|(u_{n},v_{n})\|\nonumber\\
& =&
\frac{(\alpha+\beta-2)\rho_{0}}{2(\alpha+\beta)}\|(u_{n},v_{n})\|^{2}
-\left(\frac{1}{q}-\frac{1}{\alpha+\beta}\right)(\lambda+\mu)S_{\alpha+\beta}^{q}
\max\left\{\|a\|_{\frac{\alpha+\beta}{\alpha+\beta-q}},\|c\|_{\frac{\alpha+\beta}{\alpha+\beta-q}}\right\}
\|(u_{n},v_{n})\|^{q}\nonumber\\
& &
-\left(1-\frac{1}{\alpha+\beta}\right)S_{2}(\|h_{1}\|_{2}+\|h_{2}\|_{2})\|(u_{n},v_{n})\|
\end{eqnarray}
for any $n\in \mathbb{N}$.
Hence, combining with (\ref{3.4.2}) and (\ref{3.4.3}), we mention that
\begin{eqnarray}\label{3.4.4}
\frac{(\alpha+\beta-2)\rho_{0}}{2(\alpha+\beta)}\|(u_{n},v_{n})\|^{2}
& <&
c+1
+\left(\frac{1}{q}-\frac{1}{\alpha+\beta}\right)(\lambda+\mu)S_{\alpha+\beta}^{q}
\max\left\{\|a\|_{\frac{\alpha+\beta}{\alpha+\beta-q}},\|c\|_{\frac{\alpha+\beta}{\alpha+\beta-q}}\right\}
\|(u_{n},v_{n})\|^{q}\nonumber\\
& &
+\left(1-\frac{1}{\alpha+\beta}\right)S_{2}(\|h_{1}\|_{2}+\|h_{2}\|_{2})\|(u_{n},v_{n})\|
\end{eqnarray}
for any $n>n_{0}$ and $\{(u_{n},v_{n})\}$ is bounded in $W$ since $1<q<2$.
Thus, $\|\nabla u_{n}\|_{2}$ and $\|\nabla v_{n}\|_{2}$ are bounded for all $n\in \mathbb{N}$.
Then, there is a subsequence, still denoted by $\{u_{n}\}$, such that $u_{n}\rightharpoonup u^{\ast}$ for some $u^{\ast}\in H_{0}^{1}(\Omega)$ as $n\rightarrow \infty$,
and a subsequence of $\{v_{n}\}$, still denoted by $\{v_{n}\}$, such that $v_{n}\rightharpoonup v^{\ast}$ for some $v^{\ast}\in H_{0}^{1}(\Omega)$ as $n\rightarrow \infty$.
By proposition 2.1, we know that
\begin{eqnarray}\label{3.4.5}
u_{n}\rightarrow u^{\ast} \;\mbox{in} \;L^{i_{1}}(\Omega),
\;\;
v_{n}\rightarrow v^{\ast} \;\mbox{in} \;L^{i_{2}}(\Omega),
\;\mbox{as}\; n\rightarrow \infty,
\;\mbox{for}\;1<i_{1},i_{2}<2^{\ast}.
\end{eqnarray}
Furthermore, there exist two positive constants $M_{1}$ and $M_{2}$ such that
\begin{eqnarray}\label{3.4.6}
\|u_{n}\|_{\alpha+\beta}\leq M_{1}
\;\mbox{and}\;
\|v_{n}\|_{\alpha+\beta}\leq M_{2},
\;\mbox{for all}\;n\in \mathbb{N}.
\end{eqnarray}
\par
Next, we prove that
\begin{eqnarray}\label{3.4.21}
u_{n}\rightarrow u^{\ast} \;\mbox{in} \;H_{0}^{1}(\Omega)\;\mbox{as} \; n\rightarrow\infty.
\end{eqnarray}
Let
\begin{eqnarray*}
\Gamma_{1}(u)
=\int_{\Omega}\Phi_{1}\left(\frac{u^{2}+|\nabla u|^{2}}{2}\right)dx,\;\;\forall u\in H_{0}^{1}(\Omega).
\end{eqnarray*}
It follows that
\begin{eqnarray}\label{3.4.12}
\langle \Gamma'_{1}(u), \varphi_{1}\rangle
=
\int_{\Omega}\phi_{1}\left(\frac{u^{2}+|\nabla u|^{2}}{2}\right)(u\varphi_{1}+\nabla u\cdot \nabla \varphi_{1} )dx,
\;\;\forall \varphi_{1}\in H_{0}^{1}(\Omega).
\end{eqnarray}
Then by (\ref{2.1.2}), we have
\begin{eqnarray}\label{3.4.7}
\langle J'_{\lambda,\mu}(u_{n},v_{n}),(u_{n}-u^{\ast},0)\rangle
&= &
\langle \Gamma'_{1}(u_{n}),u_{n}-u^{\ast}\rangle
-\lambda\int_{\Omega} a(x)|u_{n}|^{q-2}u_{n}(u_{n}-u^{\ast})dx
-\int_{\Omega}h_{1}(u_{n}-u^{\ast})dx\nonumber\\
& &
-\frac{\alpha}{\alpha+\beta}\int_{\Omega}b(x)|u_{n}|^{\alpha-2}u_{n}|v_{n}|^{\beta}(u_{n}-u^{\ast})dx.
\end{eqnarray}
Clearly, (\ref{3.4.1}) implies that
\begin{eqnarray}\label{3.4.8}
\langle J'_{\lambda,\mu}(u_{n},v_{n}),(u_{n}-u^{\ast},0)\rangle\rightarrow 0.
\end{eqnarray}
By $(B)$, H\"older inequality, $1<\alpha+\beta<2^{\ast}$, (\ref{3.4.5}) and (\ref{3.4.6}), we have
\begin{eqnarray}\label{3.4.9}
& &
\left|\frac{\alpha}{\alpha+\beta}\int_{\Omega}b(x)|u_{n}|^{\alpha-2}u_{n}|v_{n}|^{\beta}(u_{n}-u^{\ast})dx\right|\nonumber\\
&\leq &
\frac{\alpha}{\alpha+\beta}\int_{\Omega}|b(x)||u_{n}|^{\alpha-1}|v_{n}|^{\beta}|u_{n}-u^{\ast}|dx\nonumber\\
&\leq &
\frac{\alpha}{\alpha+\beta}\|b\|_{\infty}\int_{\Omega}|u_{n}|^{\alpha-1}|v_{n}|^{\beta}|u_{n}-u^{\ast}|dx\nonumber\\
&\leq &
\frac{\alpha}{\alpha+\beta}\|b\|_{\infty}
\left(\int_{\Omega}|u_{n}|^{(\alpha-1)\cdot\frac{\alpha+\beta}{\alpha-1}}dx\right)^{\frac{\alpha-1}{\alpha+\beta}}
\left(\int_{\Omega}|v_{n}|^{\beta\cdot\frac{\alpha+\beta}{\beta}}dx\right)^{\frac{\beta}{\alpha+\beta}}
\left(\int_{\Omega}|u_{n}-u^{\ast}|^{\alpha+\beta}dx\right)^{\frac{1}{\alpha+\beta}}\nonumber\\
&= &
\frac{\alpha}{\alpha+\beta}\|b\|_{\infty}
\|u_{n}\|_{\alpha+\beta}^{\alpha-1}
\|v_{n}\|_{\alpha+\beta}^{\beta}
\|u_{n}-u^{\ast} \|_{\alpha+\beta}\nonumber\\
&\leq &
\frac{\alpha}{\alpha+\beta}\|b\|_{\infty}
M_{1}^{\alpha-1}
M_{2}^{\beta}
\|u_{n}-u^{\ast} \|_{\alpha+\beta}
\rightarrow 0
\end{eqnarray}
as $n \rightarrow \infty$.
Then let
\begin{eqnarray*}
\delta_{1}=\frac{(\alpha+\beta)\left(\alpha+\beta+\tau_{1}(\alpha+\beta-q)\right)}
{\tau_{1}(\alpha+\beta-q)\left(\alpha+\beta+1-q\right)+\alpha+\beta}.
\end{eqnarray*}
It follows from $1<q<2<\alpha+\beta<2^{\ast}$ that $\delta_{1}>1$ and
\begin{eqnarray*}
2^{\ast}-\delta_{1}
&= & 2^{\ast}-\frac{(\alpha+\beta)\left(\alpha+\beta+\tau_{1}(\alpha+\beta-q)\right)}
{\tau_{1}(\alpha+\beta-q)\left(\alpha+\beta+1-q\right)+\alpha+\beta}\\
&= & \frac{\tau_{1}(\alpha+\beta-q)[2^{\ast}(\alpha+\beta+1-q)-(\alpha+\beta)]+[2^{\ast}-(\alpha+\beta)](\alpha+\beta)}
{\tau_{1}(\alpha+\beta-q)\left(\alpha+\beta+1-q\right)+\alpha+\beta}
>0.
\end{eqnarray*}
By the estimate just above, $(A)$, H\"older inequality, (\ref{3.4.5}) and (\ref{3.4.6}), we get
\begin{eqnarray}\label{3.4.10}
& &
\left|\int_{\Omega} a(x)|u_{n}|^{q-2}u_{n}(u_{n}-u^{\ast})dx\right|\nonumber\\
&\leq &
\int_{\Omega} |a(x)||u_{n}|^{q-1}|u_{n}-u^{\ast}|dx\nonumber\\
&\leq &
\left(\int_{\Omega}|a(x)|^{\frac{\alpha+\beta}{\alpha+\beta-q}+\tau_{1}}dx\right)
^{\frac{1}{\frac{\alpha+\beta}{\alpha+\beta-q}+\tau_{1}}}
\left(\int_{\Omega}|u_{n}|^{(q-1)\cdot\frac{\alpha+\beta}{q-1}}dx\right)^{\frac{q-1}{\alpha+\beta}}
\left(\int_{\Omega}|u_{n}-u^{\ast}|^{\delta}dx\right)^{\frac{1}{\delta_{1}}}\nonumber\\
&= &
\|a\|_{\frac{\alpha+\beta}{\alpha+\beta-q}+\tau_{1}}
\|u_{n}\|_{\alpha+\beta}^{q-1}
\|u_{n}-u^{\ast} \|_{\delta_{1}}\nonumber\\
&\leq &
\|a\|_{\frac{\alpha+\beta}{\alpha+\beta-q}+\tau_{1}}
M_{1}^{q-1}
\|u_{n}-u^{\ast} \|_{\delta_{1}}
\rightarrow 0
\end{eqnarray}
as $n \rightarrow \infty$.
Moreover, using $(H)$, H\"older inequality and (\ref{3.4.5}), we mention that
\begin{eqnarray}\label{3.4.11}
\left|\int_{\Omega}h_{1}(u_{n}-u^{\ast})dx\right|
\leq
\int_{\Omega} |h_{1}||u_{n}-u^{\ast}|dx
\leq
\|h_{1}\|_{2}
\|u_{n}-u^{\ast} \|_{2}
\rightarrow 0
\end{eqnarray}
as $n \rightarrow \infty$.
Hence, combing with (\ref{3.4.7}), (\ref{3.4.8}), (\ref{3.4.9}), (\ref{3.4.10}) and (\ref{3.4.11}),
we have
\begin{eqnarray}\label{3.4.20}
\langle \Gamma'_{1}(u_{n}),u_{n}-u^{\ast}\rangle\rightarrow 0.
\end{eqnarray}
Then combing with $(\phi_1)$-$(\phi_3)$, and applying \cite[Lemma 3.5]{Jeanjean2022}, we have
\begin{eqnarray}\label{3.4.13}
\nabla u_{n}(x)\rightarrow \nabla u^{\ast}(x) \;\mbox{a.e.}\; \;\mbox{in}\;\Omega.
 \end{eqnarray}
Next, we prove that
\begin{eqnarray}\label{3.4.14}
\limsup_{n\rightarrow\infty} \Gamma_{1}(u_{n})\leq \Gamma_{1}(u^{\ast}).
\end{eqnarray}
The proof is similar to Lemma 5.5 in \cite{Jeanjean2022}.
For the purpose of completeness, we also present the proofs here.
For $u\in H_{0}^{1}(\Omega)$ and $w\in H_{0}^{1}(\Omega)$, let $e_{1}=(u,\nabla u)$ and $e_{2}=(w,\nabla w)$ so that
\begin{eqnarray*}
e_{1}\cdot e_{2}=uw+\nabla u\cdot\nabla w.
\end{eqnarray*}
Setting $g_{1}(s)=\Phi_{1}(s^{2})$, $s\in [0,+\infty)$.
It follows that $g'_{1}(s)=2s\phi_{1}(s^{2})$.
Applying the mean value theorem, there is $\xi\in\mathbb{R}$ with
\begin{eqnarray}\label{3.4.16}
\xi\in \left(\min\left\{\frac{|e_{1}|}{\sqrt{2}},\frac{|e_{2}|}{\sqrt{2}}\right\},
\max\left\{\frac{|e_{1}|}{\sqrt{2}},\frac{|e_{2}|}{\sqrt{2}}\right\}\right)
\end{eqnarray}
such that
\begin{eqnarray*}
g_{1}\left(\frac{|e_{1}|}{\sqrt{2}}\right)-g_{1}\left(\frac{|e_{2}|}{\sqrt{2}}\right)
&= &
g'_{1}(\xi)\left(\frac{|e_{1}|-|e_{2}|}{\sqrt{2}}\right)\nonumber\\
&= &
2\xi\phi_{1}(\xi^{2})\left(\frac{|e_{1}|-|e_{2}|}{\sqrt{2}}\right)\nonumber\\
&= &
\sqrt{2}\xi\phi_{1}(\xi^{2})\left(|e_{1}|-|e_{2}|\right)\nonumber\\
&= &
\sqrt{2}\xi\phi_{1}\left(\frac{(\sqrt{2}\xi)^{2}}{2}\right)\left(|e_{1}|-|e_{2}|\right)
\end{eqnarray*}
Note that the condition $(\phi_2)$ is equivalent to requiring that the function $t\mapsto t\phi_{i}\left(\frac{t^{2}}{2}\right)$ is strictly increasing on $[0,+\infty)$.
Then, combing with (\ref{3.4.16}), we have
\begin{eqnarray*}
\sqrt{2}\xi\phi_{1}\left(\frac{(\sqrt{2}\xi)^{2}}{2}\right)\left(|e_{1}|-|e_{2}|\right)
\geq
|e_{2}|\phi_{1}\left(\frac{|e_{2}|^{2}}{2}\right)\left(|e_{1}|-|e_{2}|\right).
\end{eqnarray*}
Therefore, by the estimates above, we obtain that
\begin{eqnarray}\label{3.4.17}
\Gamma_{1}(u)-\Gamma_{1}(w)
&= &
\int_{\Omega}\left(g_{1}\left(\frac{|e_{1}|}{\sqrt{2}}\right)
-g_{1}\left(\frac{|e_{2}|}{\sqrt{2}}\right)\right)dx\nonumber\\
&\geq &
\int_{\Omega}\phi_{1}\left(\frac{|e_{2}|^{2}}{2}\right)|e_{2}|\left(|e_{1}|-|e_{2}|\right)dx.
\end{eqnarray}
In addition, it follows from (\ref{3.4.12}) that
\begin{eqnarray}\label{3.4.18}
\langle \Gamma'_{1}(w), u-w\rangle
&= &
\int_{\Omega}\phi_{1}\left(\frac{|e_{2}|^{2}}{2}\right)(w(u-w)+\nabla w\cdot \nabla (u-w))dx\nonumber\\
&= &
\int_{\Omega}\phi_{1}\left(\frac{|e_{2}|^{2}}{2}\right)e_{2}\cdot (e_{1}-e_{2})dx.
\end{eqnarray}
Thus, from (\ref{3.4.17}) and (\ref{3.4.18}), we conclude that
\begin{eqnarray}\label{3.4.19}
& &\Gamma_{1}(u)-\Gamma_{1}(w)-\langle \Gamma'_{1}(w), u-w\rangle\nonumber\\
&\geq&
\int_{\Omega}\phi_{1}\left(\frac{|e_{2}|^{2}}{2}\right)|e_{2}|\left(|e_{1}|-|e_{2}|\right)dx
-\int_{\Omega}\phi_{1}\left(\frac{|e_{2}|^{2}}{2}\right)e_{2}\cdot (e_{1}-e_{2})dx\nonumber\\
&= &
\int_{\Omega}\phi_{1}\left(\frac{|e_{2}|^{2}}{2}\right)
\left(|e_{2}||e_{1}|-e_{2}\cdot e_{1}\right)dx
\geq0
\end{eqnarray}
for all $u\in H_{0}^{1}(\Omega)$ and $w\in H_{0}^{1}(\Omega)$.
Then, let $u=u^{\ast}$ and $w=u_{n}$ in (\ref{3.4.19}), we have
\begin{eqnarray*}
\Gamma_{1}(u_{n})
\leq
\Gamma_{1}(u^{\ast})-\langle \Gamma'_{1}(u_{n}), u^{\ast}-u_{n}\rangle.
\end{eqnarray*}
Taking the limits in estimate just above, by (\ref{3.4.20}), we have (\ref{3.4.14}) holds.
Now, by $(\phi_1)$-$(\phi_3)$, (\ref{3.4.13}), (\ref{3.4.14}) and \cite[Lemma 5.4, p.18]{Jeanjean2022}, we have (\ref{3.4.21}) holds.
\par
Next, we prove that
\begin{eqnarray}\label{3.4.22}
v_{n}\rightarrow v^{\ast} \;\mbox{in} \;H_{0}^{1}(\Omega)
\;\mbox{as}\; n\rightarrow\infty.
\end{eqnarray}
Let
\begin{eqnarray*}
\Gamma_{2}(v)
=\int_{\Omega}\Phi_{2}\left(\frac{v^{2}+|\nabla v|^{2}}{2}\right)dx,\;\;\forall v\in H_{0}^{1}(\Omega).
\end{eqnarray*}
It follows that
\begin{eqnarray}\label{3.4.23}
\langle \Gamma'_{2}(v), \varphi_{2}\rangle
=
\int_{\Omega}\phi_{2}\left(\frac{v^{2}+|\nabla v|^{2}}{2}\right)(v\varphi_{2}+\nabla v\cdot \nabla \varphi_{2} )dx,
\;\;\forall \varphi_{2}\in H_{0}^{1}(\Omega).
\end{eqnarray}
Then by (\ref{2.1.2}), we have
\begin{eqnarray}\label{3.4.24}
\langle J'_{\lambda,\mu}(u_{n},v_{n}),(0,v_{n}-v^{\ast})\rangle
&= &
\langle \Gamma'_{2}(v_{n}),v_{n}-v^{\ast}\rangle
-\mu\int_{\Omega} c(x)|v_{n}|^{q-2}v_{n}(v_{n}-v^{\ast})dx
-\int_{\Omega}h_{2}(v_{n}-v^{\ast})dx\nonumber\\
& &
-\frac{\beta}{\alpha+\beta}\int_{\Omega}b(x)|v_{n}|^{\beta-2}v_{n}|u_{n}|^{\alpha}(v_{n}-v^{\ast})dx.
\end{eqnarray}
Clearly, (\ref{3.4.1}) implies that
\begin{eqnarray}\label{3.4.25}
\langle J'_{\lambda,\mu}(u_{n},v_{n}),(0,v_{n}-v^{\ast})\rangle\rightarrow 0.
\end{eqnarray}
By $(B)$, H\"older inequality, $1<\alpha+\beta<2^{\ast}$, (\ref{3.4.5}) and (\ref{3.4.6}), we have
\begin{eqnarray}\label{3.4.26}
& &
\left|\int_{\Omega}b(x)|v_{n}|^{\beta-2}v_{n}|u_{n}|^{\alpha}(v_{n}-v^{\ast})dx\right|\nonumber\\
&\leq &
\int_{\Omega}|b(x)||v_{n}|^{\beta-1}|u_{n}|^{\alpha}|v_{n}-v^{\ast}|dx\nonumber\\
&\leq &
\|b\|_{\infty}\int_{\Omega}|v_{n}|^{\beta-1}|u_{n}|^{\alpha}|v_{n}-v^{\ast}|dx\nonumber\\
&\leq &
\|b\|_{\infty}
\left(\int_{\Omega}|v_{n}|^{(\beta-1)\cdot\frac{\alpha+\beta}{\beta-1}}dx\right)^{\frac{\beta-1}{\alpha+\beta}}
\left(\int_{\Omega}|u_{n}|^{\alpha\cdot\frac{\alpha+\beta}{\alpha}}dx\right)^{\frac{\alpha}{\alpha+\beta}}
\left(\int_{\Omega}|v_{n}-v^{\ast}|^{\alpha+\beta}dx\right)^{\frac{1}{\alpha+\beta}}\nonumber\\
&= &
\|b\|_{\infty}
\|v_{n}\|_{\alpha+\beta}^{\beta-1}
\|u_{n}\|_{\alpha+\beta}^{\alpha}
\|v_{n}-v^{\ast} \|_{\alpha+\beta}\nonumber\\
&\leq &
\|b\|_{\infty}
M_{2}^{\beta-1}
M_{1}^{\alpha}
\|v_{n}-v^{\ast} \|_{\alpha+\beta}
\rightarrow 0
\end{eqnarray}
as $n \rightarrow \infty$.
Then let
\begin{eqnarray*}
\delta_{2}=\frac{(\alpha+\beta)\left(\alpha+\beta+\tau_{2}(\alpha+\beta-q)\right)}
{\tau_{2}(\alpha+\beta-q)\left(\alpha+\beta+1-q\right)+\alpha+\beta}.
\end{eqnarray*}
It follows from $1<q<2<\alpha+\beta<2^{\ast}$ that $\delta_{2}>1$ and
\begin{eqnarray*}
2^{\ast}-\delta_{2}
&= & 2^{\ast}-\frac{(\alpha+\beta)\left(\alpha+\beta+\tau_{2}(\alpha+\beta-q)\right)}
{\tau_{2}(\alpha+\beta-q)\left(\alpha+\beta+1-q\right)+\alpha+\beta}\\
&= & \frac{\tau_{2}(\alpha+\beta-q)[2^{\ast}(\alpha+\beta+1-q)-(\alpha+\beta)]+[2^{\ast}-(\alpha+\beta)](\alpha+\beta)}
{\tau_{2}(\alpha+\beta-q)\left(\alpha+\beta+1-q\right)+\alpha+\beta}
>0.
\end{eqnarray*}
By the estimate just above, $(A)$, H\"older inequality, (\ref{3.4.5}) and (\ref{3.4.6}), we get
\begin{eqnarray}\label{3.4.27}
& &
\left|\int_{\Omega} c(x)|v_{n}|^{q-2}v_{n}(v_{n}-v^{\ast})dx\right|\nonumber\\
&\leq &
\int_{\Omega} |c(x)||v_{n}|^{q-1}|v_{n}-v^{\ast}|dx\nonumber\\
&\leq &
\left(\int_{\Omega}|c(x)|^{\frac{\alpha+\beta}{\alpha+\beta-q}+\tau_{2}}dx\right)
^{\frac{1}{\frac{\alpha+\beta}{\alpha+\beta-q}+\tau_{2}}}
\left(\int_{\Omega}|v_{n}|^{(q-1)\cdot\frac{\alpha+\beta}{q-1}}dx\right)^{\frac{q-1}{\alpha+\beta}}
\left(\int_{\Omega}|v_{n}-v^{\ast}|^{\delta_{2}}dx\right)^{\frac{1}{\delta_{2}}}\nonumber\\
&= &
\|c\|_{\frac{\alpha+\beta}{\alpha+\beta-q}+\tau_{2}}
\|v_{n}\|_{\alpha+\beta}^{q-1}
\|v_{n}-v^{\ast} \|_{\delta_{2}}\nonumber\\
&\leq &
\|c\|_{\frac{\alpha+\beta}{\alpha+\beta-q}+\tau_{2}}
M_{2}^{q-1}
\|v_{n}-v^{\ast} \|_{\delta_{2}}
\rightarrow 0
\end{eqnarray}
as $n \rightarrow \infty$.
Moreover, using $(H)$, H\"older inequality and (\ref{3.4.5}), we mention that
\begin{eqnarray}\label{3.4.28}
\left|\int_{\Omega}h_{2}(v_{n}-v^{\ast})dx\right|
\leq
\int_{\Omega} |h_{2}||v_{n}-v^{\ast}|dx
\leq
\|h_{2}\|_{2}
\|v_{n}-v^{\ast} \|_{2}
\rightarrow 0
\end{eqnarray}
as $n \rightarrow \infty$.
Hence, by (\ref{3.4.24})-(\ref{3.4.28}),
we have
\begin{eqnarray}\label{3.4.29}
\langle \Gamma'_{2}(v_{n}),v_{n}-v^{\ast}\rangle\rightarrow 0.
\end{eqnarray}
Then, combing with $(\phi_1)$-$(\phi_3)$, and applying \cite[Lemma 3.5, p.11]{Jeanjean2022}, we obtain that
\begin{eqnarray}\label{3.4.30}
\nabla v_{n}(x)\rightarrow \nabla v^{\ast}(x) \;\mbox{a.e.}\; \;\mbox{in}\;\Omega.
 \end{eqnarray}
Moreover, similarly the proof of (\ref{3.4.14}), by (\ref{3.4.29}) and $(\phi_2)$, we also have
\begin{eqnarray}\label{3.4.31}
\limsup_{n\rightarrow\infty} \Gamma_{2}(v_{n})\leq \Gamma_{2}(v^{\ast}).
\end{eqnarray}
Now, by $(\phi_1)$-$(\phi_3)$, (\ref{3.4.30}), (\ref{3.4.31}) and \cite[Lemma 5.4, p.18]{Jeanjean2022}, we have (\ref{3.4.22}) holds.
The proof for this lemma now is completed.
\qed
 \vskip2mm
 \noindent
{\bf Lemma 3.5.} Assume that $(\phi_1)$, $(A)$, $(B)$ and $(H)$ hold.
Then for each $\lambda+\mu\in (0,\Lambda_{0})$ and $h_{1},h_{2}\in L^{2}(\Omega)$ with
$\|h_{1}\|_{2}^{2}+\|h_{2}\|_{2}^{2}\in (0, m_{\lambda,\mu}]$, we have
\begin{eqnarray*}
-\infty<\inf\{J_{\lambda,\mu}(u,v):(u,v)\in \bar{B}_{\rho_{\lambda,\mu}}\}<0,
\end{eqnarray*}
where $\bar{B}_{\rho_{\lambda,\mu}}=\{(u,v)\in W:\|(u,v)\|\leq\rho_{\lambda,\mu}\}$ and $\rho_{\lambda,\mu}$ is given in Lemma 3.2.
\vskip2mm
  \noindent
{\bf Proof.} Since $h_{1},h_{2}\in L^{2}(\Omega)$ and $h_{1},h_{2}\not\equiv0$, we can choose a function $(\tilde{\varphi}_{1},\tilde{\varphi}_{2})\in C_{0}^{\infty}(\Omega)\times C_{0}^{\infty}(\Omega)$ such that
\begin{eqnarray}\label{3.5.1}
\int_{\Omega}h_{1}(x)\tilde{\varphi}_{1}(x)dx>0
\;\mbox{and}\;
\int_{\Omega}h_{2}(x)\tilde{\varphi}_{2}(x)dx>0.
\end{eqnarray}
Then by (\ref{2.1.1}) and (\ref{3.1.2}), we have
\begin{eqnarray*}
J_{\lambda,\mu}(t\tilde{\varphi}_{1},t\tilde{\varphi}_{2})
& \leq&
t^{2}\frac{\rho_{1}}{2}\|(\tilde{\varphi}_{1},\tilde{\varphi}_{2})\|^{2}
-
t^{q}\frac{1}{q}\int_{\Omega}(\lambda a(x)|\tilde{\varphi}_{1}|^{q}+\mu c(x)|\tilde{\varphi}_{2}|^{q})dx
-
t^{\alpha+\beta}\frac{1}{\alpha+\beta}\int_{\Omega}b(x)|\tilde{\varphi}_{1}|^{\alpha}|\tilde{\varphi}_{2}|^{\beta}dx
\nonumber\\
& &
-t\int_{\Omega}(h_{1}\tilde{\varphi}_{1}+h_{2}\tilde{\varphi}_{2})dx.
\end{eqnarray*}
Note that $1<q<2<\alpha+\beta<2^{\ast}$.
Using (\ref{3.5.1}), there exists a sufficiently small $\tilde{t}>0$ such that $\|(\tilde{t}\varphi_{1},\tilde{t}\varphi_{2})\|<\rho_{\lambda,\mu}$ and $J_{\lambda,\mu}(\tilde{t}\varphi_{1},\tilde{t}\varphi_{2})<0$.
Hence,
\begin{eqnarray*}
\inf\{J_{\lambda,\mu}(u,v):(u,v)\in \bar{B}_{\rho_{\lambda,\mu}}\}
\leq
J_{\lambda,\mu}(\tilde{t}\varphi_{1},\tilde{t}\varphi_{2})
<0.
\end{eqnarray*}
In addition, it is easy to see that (\ref{3.2.3}) still holds for all $(u,v)\in \bar{B}_{\rho_{\lambda,\mu}}$.
Then by $\epsilon_{1}^{2}+\epsilon_{2}^{2}<\rho_{0}/S_{2}$ and $\|h_{1}\|_{2}^{2}+\|h_{2}\|_{2}^{2}\leq m_{\lambda,\mu}$, we have
\begin{eqnarray}\label{3.5.2}
J_{\lambda,\mu}(u,v)
& \geq&
\left(\frac{\rho_{0}}{2}-\frac{S_{2}}{2}\left(\epsilon_{1}^{2}+\epsilon_{2}^{2}\right)\right)\|(u,v)\|^{2}
-
(\lambda+\mu)\frac{1}{q}S_{\alpha+\beta}^{q}
\max\left\{\|a\|_{\frac{\alpha+\beta}{\alpha+\beta-q}},\|c\|_{\frac{\alpha+\beta}{\alpha+\beta-q}}\right\}
\|(u,v)\|^{q}\nonumber\\
& &
-\frac{1}{\alpha+\beta}\|b\|_{\infty}S_{\alpha+\beta}^{\alpha+\beta}\|(u,v)\|^{\alpha+\beta}
-\max\left\{\frac{S_{2}}{2\epsilon_{1}^{2}},\frac{S_{2}}{2\epsilon_{2}^{2}}\right\}
\left(\|h_{1}\|_{2}^{2}+\|h_{2}\|_{2}^{2}\right)\nonumber\\
& \geq&
-
(\lambda+\mu)\frac{1}{q}S_{\alpha+\beta}^{q}
\max\left\{\|a\|_{\frac{\alpha+\beta}{\alpha+\beta-q}},\|c\|_{\frac{\alpha+\beta}{\alpha+\beta-q}}\right\}
\|(u,v)\|^{q}
-\frac{1}{\alpha+\beta}\|b\|_{\infty}S_{\alpha+\beta}^{\alpha+\beta}\|(u,v)\|^{\alpha+\beta}\nonumber\\
& &-\max\left\{\frac{S_{2}}{2\epsilon_{1}^{2}},\frac{S_{2}}{2\epsilon_{2}^{2}}\right\}
m_{\lambda,\mu}\nonumber\\
& \geq&
-
(\lambda+\mu)\frac{1}{q}S_{\alpha+\beta}^{q}
\max\left\{\|a\|_{\frac{\alpha+\beta}{\alpha+\beta-q}},\|c\|_{\frac{\alpha+\beta}{\alpha+\beta-q}}\right\}
\rho_{\lambda,\mu}^{q}
-\frac{1}{\alpha+\beta}\|b\|_{\infty}S_{\alpha+\beta}^{\alpha+\beta}\rho_{\lambda,\mu}^{\alpha+\beta}\nonumber\\
& &-\max\left\{\frac{1}{2\epsilon_{S_{2}}^{2}},\frac{S_{2}}{2\epsilon_{2}^{2}}\right\}
m_{\lambda,\mu},
\end{eqnarray}
which implies that $J_{\lambda,\mu}$ is bounded from blew in $\bar{B}_{\rho_{\lambda,\mu}}$ for each $\lambda+\mu\in (0,\Lambda_{0})$ and $h_{1},h_{2}\in L^{2}(\Omega)$ with
$\|h_{1}\|_{2}^{2}+\|h_{2}\|_{2}^{2}\in (0, m_{\lambda,\mu}]$.
So
$\inf\{J_{\lambda,\mu}(u,v):(u,v)\in \bar{B}_{\rho_{\lambda,\mu}}\}>-\infty$.
The proof for this lemma now complete.
\qed

\vskip2mm
 \noindent
{\bf The proof of Theorem 1.1. } By Lemma 3.2-Lemma 3.4 and Lemma 2.3, we obtain that for each $\lambda+\mu\in (0,\Lambda_{0})$ and $h_{1},h_{2}\in L^{2}(\Omega)$ with
$\|h_{1}\|_{2}^{2}+\|h_{2}\|_{2}^{2}\in (0, m_{\lambda,\mu}]$,
$J_{\lambda,\mu}$ has a critical value $c_{\ast}\geq\alpha_{\lambda,\mu}>0$ with
 \begin{eqnarray*}
 c_{\ast}:=\inf_{\gamma\in\Gamma}\max_{t\in[0,1]}\varphi(\gamma(t)),
\end{eqnarray*}
where
\begin{eqnarray*}
\Gamma:=\{\gamma\in C([0,1],X):\gamma(0)=(0,0), \; \gamma(1)=(t_{1}\varphi_{1},t_{1}\varphi_{2})\},
\end{eqnarray*}
$\alpha_{\lambda,\mu}$ and $(t_{1}\varphi_{1},t_{1}\varphi_{2})$ were given in Lemma 3.2 and Lemma 3.3, respectively.
Hence, system (\ref{eq1}) has one nontrivial solution $(u_{0},v_{0})$ of positive energy.
Obviously, $(u_{0},v_{0})\neq (0,0)$. Otherwise, by the fact that $\Phi_{i}(0)=0$, $i=1,2$, we have $J_{\lambda,\mu}(u_{0},v_{0})=0$, which contradicts with $c_{\ast}>0$.
\par
Next, we prove that system (\ref{eq1}) has one nontrivial solution of negative energy.
By Lemma 3.2 and Lemma 3.5, we conclude that
\begin{eqnarray*}
-\infty<\inf_{\bar{B}_{\rho_{\lambda,\mu}}}J_{\lambda,\mu}
<0<\inf_{\partial B_{\rho_{\lambda,\mu}}}J_{\lambda,\mu}
\end{eqnarray*}
for each $\lambda+\mu\in (0,\Lambda_{0})$ and $h_{1},h_{2}\in L^{2}(\Omega)$ with
$\|h_{1}\|_{2}^{2}+\|h_{2}\|_{2}^{2}\in (0, m_{\lambda,\mu}]$,
where $\partial B_{\rho_{\lambda,\mu}}=\{(u,v)\in W:\|(u,v)\|=\rho_{\lambda,\mu}\}$,
$\Lambda_{0}$, $m_{\lambda,\mu}$ and $\rho_{\lambda,\mu}$ were given in Lemma 3.2.
Taking
\begin{eqnarray}\label{1.1.1}
\frac{1}{n}
\in
\left(0, \inf_{\partial B_{\rho_{\lambda,\mu}}}J_{\lambda,\mu}
-\inf_{\bar{B}_{\rho_{\lambda,\mu}}}J_{\lambda,\mu}\right),
\;\;
n\in \mathbb{N}^{+}.
\end{eqnarray}
It follows from the definition of $\inf_{\bar{B}_{\rho_{\lambda,\mu}}}J_{\lambda,\mu}$ that there exists a $(u_{n},v_{n})\in \bar{B}_{\rho_{\lambda,\mu}}$ such that
\begin{eqnarray}\label{1.1.2}
J_{\lambda,\mu}(u_{n},v_{n})
\leq
\inf_{\bar{B}_{\rho_{\lambda,\mu}}}J_{\lambda,\mu}+\frac{1}{n}.
\end{eqnarray}
As $J_{\lambda,\mu}$ is of class $C^{1}$, we know that $J_{\lambda,\mu}$ is lower semicontinuous.
Therefore, by Lemma 2.1 we obtain that
\begin{eqnarray}\label{1.1.3}
J_{\lambda,\mu}(u_{n},v_{n})
\leq
J_{\lambda,\mu}(u,v)
+\frac{1}{n}\|(u,v)-(u_{n},v_{n})\|,
\;\;
\forall (u,v)\in \bar{B}_{\rho_{\lambda,\mu}}.
\end{eqnarray}
It follows from (\ref{1.1.1}) and (\ref{1.1.2}) that
\begin{eqnarray*}
J_{\lambda,\mu}(u_{n},v_{n})
\leq
\inf_{\bar{B}_{\rho_{\lambda,\mu}}}J_{\lambda,\mu}+\frac{1}{n}
<\inf_{\partial B_{\rho_{\lambda,\mu}}}J_{\lambda,\mu}.
\end{eqnarray*}
Hence, $(u_{n},v_{n})\in B_{\rho_{\lambda,\mu}}$.
Define $M_{n}:\bar{B}_{\rho_{\lambda,\mu}}\rightarrow\mathbb{R}$ by
\begin{eqnarray*}
M_{n}(u,v)=J_{\lambda,\mu}(u,v)+\frac{1}{n}\|(u,v)-(u_{n},v_{n})\|.
\end{eqnarray*}
Clearly, by (\ref{1.1.3}), we have $M_{n}(u_{n},v_{n})=J_{\lambda,\mu}(u_{n},v_{n})\leq M_{n}(u,v)$,
which implies that $(u_{n},v_{n})\in B_{\rho_{\lambda,\mu}}$ minimizes $M_{n}$ on $\bar{B}_{\rho_{\lambda,\mu}}$.
Therefore, for all $(\varphi_{1},\varphi_{2})\in W$ with $\|(\varphi_{1},\varphi_{2})\|=1$, taking $t>0$ small enough such that $(u_{n}+t\varphi_{1},v_{n}+t\varphi_{2})\in \bar{B}_{\rho_{\lambda,\mu}}$, then
\begin{eqnarray*}
0
& \leq &
\frac{M_{n}(u_{n}+t\varphi_{1},v_{n}+t\varphi_{2})-M_{n}(u_{n},v_{n})}{t}\\
&= &
\frac{J_{\lambda,\mu}(u_{n}+t\varphi_{1},v_{n}+t\varphi_{2})
+\frac{1}{n}\|(u_{n}+t\varphi_{1},v_{n}+t\varphi_{2})-(u_{n},v_{n})\|-J_{\lambda,\mu}(u_{n},v_{n})}{t}\\
&= &
\frac{J_{\lambda,\mu}(u_{n}+t\varphi_{1},v_{n}+t\varphi_{2})-J_{\lambda,\mu}(u_{n},v_{n})}{t}
+\frac{\|(u_{n},v_{n})+t(\varphi_{1},\varphi_{2})-(u_{n},v_{n})\|}{nt}\\
&= &
\frac{J_{\lambda,\mu}(u_{n}+t\varphi_{1},v_{n}+t\varphi_{2})-J_{\lambda,\mu}(u_{n},v_{n})}{t}
+\frac{1}{n},
\end{eqnarray*}
which implies that
\begin{eqnarray*}
\frac{J_{\lambda,\mu}(u_{n}+t\varphi_{1},v_{n}+t\varphi_{2})-J_{\lambda,\mu}(u_{n},v_{n})}{t}
\geq-\frac{1}{n}.
\end{eqnarray*}
Taking the limits in estimate just above, we have
\begin{eqnarray}\label{1.1.4}
\langle J'_{\lambda,\mu}(u_{n},v_{n}),(\varphi_{1},\varphi_{2})\rangle
=
\lim_{t\rightarrow0}
\frac{J_{\lambda,\mu}(u_{n}+t\varphi_{1},v_{n}+t\varphi_{2})-J_{\lambda,\mu}(u_{n},v_{n})}{t}
\geq-\frac{1}{n}
\end{eqnarray}
Similarly, when $t<0$ and $|t|$ small enough, we conclude that
\begin{eqnarray}\label{1.1.5}
\langle J'_{\lambda,\mu}(u_{n},v_{n}),(\varphi_{1},\varphi_{2})\rangle
\leq\frac{1}{n}.
\end{eqnarray}
Combing with (\ref{1.1.4}) and (\ref{1.1.5}), we deduce that
\begin{eqnarray}\label{1.1.6}
\|J'_{\lambda,\mu}(u_{n},v_{n})\|_{(H^{1}_{0}(\Omega))^{'}}\leq\frac{1}{n}.
\end{eqnarray}
Now, passing to the limit in (\ref{1.1.2}) and (\ref{1.1.6}), we mention that
\begin{eqnarray*}
\lim_{n\rightarrow\infty}J_{\lambda,\mu}(u_{n},v_{n})=\inf_{\bar{B}_{\rho_{\lambda,\mu}}}J_{\lambda,\mu}<0
\;\;\mbox{and} \;
\lim_{n\rightarrow\infty}J'_{\lambda,\mu}(u_{n},v_{n})=0.
\end{eqnarray*}
Thus, $\{(u_{n},v_{n})\}\subset \bar{B}_{\rho_{\lambda,\mu}}$ is a Palais-Smale sequence of $J_{\lambda,\mu}$.
By Lemma 3.4, $\{(u_{n},v_{n})\}$ has a strongly convergent subsequence
$\{(u_{n_{i}},v_{n_{i}})\}\subset \bar{B}_{\rho_{\lambda,\mu}}$, and
$(u_{n_{i}},v_{n_{i}})\rightarrow (\tilde{u}_{0},\tilde{v}_{0})\in\bar{B}_{\rho_{\lambda,\mu}}$ as $n_{i}\rightarrow\infty$.
As a consequence,
\begin{eqnarray*}
J_{\lambda,\mu}(\tilde{u}_{0},\tilde{v}_{0})=\inf_{\bar{B}_{\rho_{\lambda,\mu}}}J_{\lambda,\mu}<0
\;\;\mbox{and} \;
J'_{\lambda,\mu}(\tilde{u}_{0},\tilde{v}_{0})=0,
\end{eqnarray*}
which imply that problem (\ref{eq1}) has one solution $(\tilde{u}_{0},\tilde{v}_{0})$ of negative energy.
Obviously, $(\tilde{u}_{0},\tilde{v}_{0})\neq(0,0)$.
Otherwise, by the fact that $\Phi_{i}(0)=0$, $i=1,2$, we have $J_{\lambda,\mu}(\tilde{u}_{0},\tilde{v}_{0})=0$, which contradicts with $\inf_{\bar{B}_{\rho_{\lambda,\mu}}}J_{\lambda,\mu}<0$.
\par
Finally, we investigate the non-semi-trivial solution of problem (\ref{eq1}).
If the nontrivial solution $(u,v)=(u,0)$, then by (\ref{2.1.2}), there holds
\begin{eqnarray}\label{rr1}
\int_{\Omega}\phi_{1}\left(\frac{u^{2}+|\nabla u|^{2}}{2}\right)(u^{2}+|\nabla u|^{2})dx
-\lambda\int_{\Omega}a(x)|u|^{q}dx
-\int_{\Omega}h_{1}udx=0.
\end{eqnarray}
If $\lambda+\mu\in (0,\Lambda_{2})$ and $h_{1},h_{2}\in L^{2}(\Omega)$ with $\|h_{1}\|_{2}^{2}+\|h_{2}\|_{2}^{2}\in (0, m^{\ast}_{\lambda,\mu}]$,
by $(\phi_1)$, (\ref{3.1.3}), (\ref{3.2.1}) and $\epsilon_{1}^{2}+\epsilon_{2}^{2}<\rho_{0}/S_{2}$, we have
\begin{eqnarray*}
& &
\int_{\Omega}\phi_{1}\left(\frac{u^{2}+|\nabla u|^{2}}{2}\right)(u^{2}+|\nabla u|^{2})dx
-\lambda\int_{\Omega}a(x)|u|^{q}dx
-\int_{\Omega}h_{1}udx\\
&> &
\rho_{0}\|\nabla u\|_{2}^{2}
-\lambda S_{\alpha+\beta}^{q}\|a\|_{\frac{\alpha+\beta}{\alpha+\beta-q}}\|\nabla u\|_{2}^{q}
-S_{2}\left(\frac{\epsilon_{1}^{2}}{2}\|\nabla u\|_{2}^{2}+\frac{1}{2\epsilon_{1}^{2}}\|h_{1}\|_{2}^{2}\right)\\
&= &
\frac{2\rho_{0}-S_{2}\epsilon_{1}^{2}}{2}\|\nabla u\|_{2}^{2}
-\lambda S_{\alpha+\beta}^{q}\|a\|_{\frac{\alpha+\beta}{\alpha+\beta-q}}\|\nabla u\|_{2}^{q}
-\frac{S_{2}}{2\epsilon_{1}^{2}}\|h_{1}\|_{2}^{2}>0,
\end{eqnarray*}
which contradicts with (\ref{rr1}).
Hence, the nontrivial solution $(u,v)\neq(u,0)$.
Similarly, we have the nontrivial solution $(u,v)\neq(0,v)$.
As a consequence, both $(u_{0},v_{0})$ and $(\tilde{u}_{0},\tilde{v}_{0})$ are non-semi-trivial solutions of system (\ref{eq1}).
The proof is completed.
\qed
\vskip2mm
 \noindent
{\bf The proof of Theorem 1.2. } Note that $(\phi_{4})$ can be deleted if we strengthen $(\phi_{1})$ to $(\phi_1)'$.
Observe the proof of Theorem 1.1, we find that $(\phi_{4})$ is only used to guarantee the boundedness of Palais-Smale sequences. Hence,
here we just need to prove that any Palais-Smale sequence $\{(u_{n},v_{n})\}$ of $J_{\lambda,\mu}$ is bounded in $W$.
Let $c\in\mathbb{R}$ and $\{(u_{n},v_{n})\}$ be a $(PS)_{c}$ sequence of $J_{\lambda,\mu}$ in $W$, that is,
\begin{eqnarray}\label{5.4.1}
\lim_{n\rightarrow\infty}J_{\lambda,\mu}(u_{n},v_{n})=c,\;\;\lim_{n\rightarrow\infty}J'_{\lambda,\mu}(u_{n},v_{n})=0.
\end{eqnarray}
We first claim that $\{(u_{n},v_{n})\}$ is bounded in $W$.
Clearly, by (\ref{5.4.1}), we have
\begin{eqnarray*}
J_{\lambda,\mu}(u_{n},v_{n})-\frac{1}{\alpha+\beta}\langle J'_{\lambda,\mu}(u_{n},v_{n}),(u_{n},v_{n}) \rangle
\rightarrow
c
<
c+1,\;\;\mbox{as}\;\; n\rightarrow\infty.
\end{eqnarray*}
Thus, there exists $n_{0}>0$ large enough such that
\begin{eqnarray}\label{5.4.2}
J_{\lambda,\mu}(u_{n},v_{n})-\frac{1}{\alpha+\beta}\langle J'_{\lambda,\mu}(u_{n},v_{n}),(u_{n},v_{n}) \rangle
<
c+1
\end{eqnarray}
for any $n>n_{0}$.
In addition, by $(\phi_1)'$, $(A)$, $(H)$, (\ref{2.1.1}), (\ref{3.1.1}) and (\ref{3.1.3}), we have
\begin{eqnarray}\label{5.4.3}
& &
J_{\lambda,\mu}(u_{n},v_{n})-\frac{1}{\alpha+\beta}\langle J'_{\lambda,\mu}(u_{n},v_{n}),(u_{n},v_{n}) \rangle\nonumber\\
&= &
\int_{\Omega}\left[\Phi_{1}\left(\frac{u_{n}^{2}+|\nabla u_{n}|^{2}}{2}\right)
-\frac{1}{\alpha+\beta}\phi_{1}\left(\frac{u_{n}^{2}+|\nabla u_{n}|^{2}}{2}\right)(u_{n}^{2}+|\nabla u_{n}|^{2})\right]dx\nonumber\\
& &
+\int_{\Omega}\left[\Phi_{2}\left(\frac{v_{n}^{2}+|\nabla v_{n}|^{2}}{2}\right)
-\frac{1}{\alpha+\beta}\phi_{2}\left(\frac{v_{n}^{2}+|\nabla v_{n}|^{2}}{2}\right)(v_{n}^{2}+|\nabla v_{n}|^{2})\right]dx\nonumber\\
& &
-\left(\frac{1}{q}-\frac{1}{\alpha+\beta}\right)\int_{\Omega}(\lambda a(x)|u_{n}|^{q}+\mu c(x)|v_{n}|^{q})dx
-\left(1-\frac{1}{\alpha+\beta}\right)\int_{\Omega}(h_{1}u_{n}+h_{2}v_{n})dx
\nonumber\\
& \geq&
\int_{\Omega}\left[\rho_0\frac{u_{n}^{2}+|\nabla u_{n}|^{2}}{2}
-\frac{\rho_1}{\alpha+\beta}(u_{n}^{2}+|\nabla u_{n}|^{2})\right]dx
+\int_{\Omega}\left[\rho_0\frac{v_{n}^{2}+|\nabla v_{n}|^{2}}{2}
-\frac{\rho_1}{\alpha+\beta}(v_{n}^{2}+|\nabla v_{n}|^{2})\right]dx\nonumber\\
& &
-\left(\frac{1}{q}-\frac{1}{\alpha+\beta}\right)(\lambda+\mu)S_{\alpha+\beta}^{q}
\max\left\{\|a\|_{\frac{\alpha+\beta}{\alpha+\beta-q}},\|c\|_{\frac{\alpha+\beta}{\alpha+\beta-q}}\right\}
\|(u_{n},v_{n})\|^{q}
-\left(1-\frac{1}{\alpha+\beta}\right)S_{2}(\|h_{1}\|_{2}+\|h_{2}\|_{2})\|(u_{n},v_{n})\|
\nonumber\\
& =&
\left(\frac{\rho_0}{2}-\frac{\rho_1}{\alpha+\beta}\right)\int_{\Omega}(u_{n}^{2}+|\nabla u_{n}|^{2}+v_{n}^{2}+|\nabla v_{n}|^{2})dx
-\left(1-\frac{1}{\alpha+\beta}\right)S_{2}(\|h_{1}\|_{2}+\|h_{2}\|_{2})\|(u_{n},v_{n})\|\nonumber\\
& &
-\left(\frac{1}{q}-\frac{1}{\alpha+\beta}\right)(\lambda+\mu)S_{\alpha+\beta}^{q}
\max\left\{\|a\|_{\frac{\alpha+\beta}{\alpha+\beta-q}},\|c\|_{\frac{\alpha+\beta}{\alpha+\beta-q}}\right\}
\|(u_{n},v_{n})\|^{q}
\nonumber\\
& \geq&
\left(\frac{\rho_0}{2}-\frac{\rho_1}{\alpha+\beta}\right)\int_{\Omega}(|\nabla u_{n}|^{2}+|\nabla v_{n}|^{2})dx
-\left(\frac{1}{q}-\frac{1}{\alpha+\beta}\right)(\lambda+\mu)S_{\alpha+\beta}^{q}
\max\left\{\|a\|_{\frac{\alpha+\beta}{\alpha+\beta-q}},\|c\|_{\frac{\alpha+\beta}{\alpha+\beta-q}}\right\}
\|(u_{n},v_{n})\|^{q}\nonumber\\
& &
-\left(1-\frac{1}{\alpha+\beta}\right)S_{2}(\|h_{1}\|_{2}+\|h_{2}\|_{2})\|(u_{n},v_{n})\|\nonumber\\
& =&
\left(\frac{\rho_0}{2}-\frac{\rho_1}{\alpha+\beta}\right)\|(u_{n},v_{n})\|^{2}
-\left(\frac{1}{q}-\frac{1}{\alpha+\beta}\right)(\lambda+\mu)S_{\alpha+\beta}^{q}
\max\left\{\|a\|_{\frac{\alpha+\beta}{\alpha+\beta-q}},\|c\|_{\frac{\alpha+\beta}{\alpha+\beta-q}}\right\}
\|(u_{n},v_{n})\|^{q}\nonumber\\
& &
-\left(1-\frac{1}{\alpha+\beta}\right)S_{2}(\|h_{1}\|_{2}+\|h_{2}\|_{2})\|(u_{n},v_{n})\|
\end{eqnarray}
for any $n\in \mathbb{N}$.
Hence, combining with (\ref{5.4.2}) and (\ref{5.4.3}), we mention that
\begin{eqnarray*}
\left(\frac{\rho_0}{2}-\frac{\rho_1}{\alpha+\beta}\right)\|(u_{n},v_{n})\|^{2}
& <&
c+1
+\left(\frac{1}{q}-\frac{1}{\alpha+\beta}\right)(\lambda+\mu)S_{\alpha+\beta}^{q}
\max\left\{\|a\|_{\frac{\alpha+\beta}{\alpha+\beta-q}},\|c\|_{\frac{\alpha+\beta}{\alpha+\beta-q}}\right\}
\|(u_{n},v_{n})\|^{q}\nonumber\\
& &
+\left(1-\frac{1}{\alpha+\beta}\right)S_{2}(\|h_{1}\|_{2}+\|h_{2}\|_{2})\|(u_{n},v_{n})\|
\end{eqnarray*}
for any $n>n_{0}$ and $\{(u_{n},v_{n})\}$ is bounded in $W$ since $1<q<2$.
\qed
\section{Examples}
\par
In this section, we present some examples to illustrate our main results.
For system (\ref{eq1}), $\phi_{i}(i=1,2)$ can be chosen from the following cases which satisfy all $(\phi_1)$-$(\phi_4)$ type conditions:
\begin{itemize}
 \item[$1.$] Let $\Phi(s)=\ln(1+s)+As$ with $A>1$.
 \item[$2.$] Let $\Phi(s)=As-\frac{1}{1-\alpha}(1+s)^{1-\alpha}$ with $\alpha>1$ and $A>1$.
 \item[$3.$] Let $\Phi(s)=As+\frac{1}{1-\alpha}(1+s)^{1-\alpha}$ with $\alpha\in[0,\frac{1}{2}]$, $A>0$
 or $\alpha\in(\frac{1}{2},1)$, $A>2\alpha-1$.
 \end{itemize}
\par
Then, we also give some cases that satisfy $(\phi_1)$-$(\phi_3)$ but do not satisfy $(\phi_4)$ type condition:
\begin{itemize}
 \item[$4.$] Let $\Phi(s)=As-\ln(1+s)$ with $A>1$.
 \item[$5.$] Let $\Phi(s)=As-\frac{1}{1-\alpha}(1+s)^{1-\alpha}$ with $0\leq\alpha<1$ and $A>1$.
 \item[$6.$] Let $\Phi(s)=As+\frac{1}{1-\alpha}(1+s)^{1-\alpha}$ with $\alpha>1$, $A>2\alpha-1$.
 \item[$7.$] Let $\Phi(s)=As+\frac{1}{1-\alpha}(1+s)^{1-\alpha}-\frac{1}{1-\beta}(1+s)^{1-\beta}$ with $0\leq\alpha\leq\frac{1}{2}<\beta<1$, $A>0$
     or $0<\beta<1<\alpha$, $A>2\alpha$
    or $1<\alpha<\beta$, $A>2\alpha$.
 \end{itemize}
 \par
Moreover, for system (\ref{eq1}), $\phi_{i}(i=1,2)$ can be chosen from the following cases which satisfy $(\phi_1)'$, $(\phi_2)$-$(\phi_4)$ type conditions:
\begin{itemize}
 \item[$8.$] Let $\Phi(s)=\ln(1+s)+As$ with $A>\frac{2}{\alpha+\beta-2}$.
 \item[$9.$] Let $\Phi(s)=As-\frac{1}{1-\alpha}(1+s)^{1-\alpha}$ with $\alpha>1$ and $A>\frac{\alpha+\beta}{\alpha+\beta-2}$.
 \item[$10.$] Let $\Phi(s)=As+\frac{1}{1-\alpha}(1+s)^{1-\alpha}$ with $\alpha\in[0,\frac{1}{2}]$, $A>\frac{2}{\alpha+\beta-2}$
 or $\alpha\in(\frac{1}{2},1)$, $A>\max\{2\alpha-1,\frac{2}{\alpha+\beta-2}\}$.
 \end{itemize}
\par
Finally, we also give some cases that satisfy $(\phi_1)'$, $(\phi_2)$ and $(\phi_3)$ but do not satisfy $(\phi_4)$ type condition:
\begin{itemize}
 \item[$11.$] Let $\Phi(s)=As-\ln(1+s)$ with $A>\frac{\alpha+\beta}{\alpha+\beta-2}$.
 \item[$12.$] Let $\Phi(s)=As-\frac{1}{1-\alpha}(1+s)^{1-\alpha}$ with $0\leq\alpha<1$ and $A>\frac{\alpha+\beta}{\alpha+\beta-2}$.
 \item[$13.$] Let $\Phi(s)=As+\frac{1}{1-\alpha}(1+s)^{1-\alpha}$ with $\alpha>1$, $A>\max\{2\alpha-1,\frac{2}{\alpha+\beta-2}\}$.
 \item[$14.$] Let $\Phi(s)=As+\frac{1}{1-\alpha}(1+s)^{1-\alpha}-\frac{1}{1-\beta}(1+s)^{1-\beta}$ with $0\leq\alpha\leq\frac{1}{2}<\beta<1$, $A>\frac{\alpha+\beta+2}{\alpha+\beta-2}$
     or $0<\beta<1<\alpha$, $A>\max\{2\alpha,\frac{\alpha+\beta+2}{\alpha+\beta-2}\}$
    or $1<\alpha<\beta$, $A>\max\{2\alpha,\frac{\alpha+\beta+2}{\alpha+\beta-2}\}$.
 \end{itemize}

\vskip2mm
 \noindent
{\bf Acknowledgements}\\
This work is supported by Yunnan Fundamental Research Projects (grant No: 202301AT070465) of China and  Xingdian Talent
Support Program for Young Talents of Yunnan Province in China.
\vskip2mm
 \noindent
 {\bf Statements and Declarations}\\
The authors state that there is no conflict of interest.\\

\renewcommand\refname{References}


\begin{thebibliography}{99}
\bibitem{Ambrosetti1994} Ambrosetti, A., Brezis, H., Cerami, G.: Combined effects of concave and convex nonlinearities in some elliptic problems. J. Funct. Anal., 122 (1994), 519-543.
\bibitem{Adams2003} Adames, R.A., Fournier, J.J.F.: Sobolev Spaces. Academic Press., 2003.
\bibitem{Afrouzi2009} Afrouzi, G.A., Rasouli, S.H.: A remark on the existence and multiplicity result for a nonlinear elliptic problem involving the $p$-Laplacian. NoDEA-Nonlinear Diff., 16 (2009), 717-730.
\bibitem{Brown2007}  Brown, K.J., Wu, T.F.: A fibering map approach to a semilinear elliptic boundary value problem. Electron. J. Differ. Eq., 69 (2007), 1-9.
\bibitem{Brown2009} Brown, K.J., Wu, T.F.: A fibering map approach to a potential operator equation and its applications. Differ. Integ. Equ., 22 (2009), 1097-1114.
\bibitem{Chen2016} Chen, W.J., Xie, J.H.: Multiple solutions to the nonhomogeneous Kirchhoff type problem involving a nonlocal operator. Differ. Equat. Appl., 3 (2016), 367-375.
\bibitem{Chen1991} Chen, Y.: TE and TM families of self-trapped beams. IEEE J. Quantum Elect., 27 (1991),
1236-1241.
\bibitem{Snyder1991} Chen, Y., Snyder, A.W.: TM-type self-guided beams with circular cross-section. Electron.
Lett., 27 (1991), 564-566.
\bibitem{Chen2013} Chen, C.S., Huang, J.C., Liu, L.H.: Multiple solutions to the nonhomogeneous $p$-Kirchhoff elliptic equation with concave-convex nonlinearities. Appl. Math. Lett., 26 (2013), 754-759.
\bibitem{Chen2011} Chen, C.Y., Kuo, Y.C., Wu, T.F.: The Nehari manifold for a Kirchhoff type problem involving sign-changing weight functions. J. Differ. Equations, 250 (2011), 1876-1908.
\bibitem{Chen2012} Chen, C.Y., Wu, T.F.: The Nehari manifold for indefinite semilinear elliptic systems involving critical exponent. Appl. Math. Comput., 218 (2012), 10817-10828.
\bibitem{Chen2022} Chen, G.F., Wu, T.F.: Multiple positive solutions for a class of Kirchhoff type equations with indefinite nonlinearities. Adv. Nonlinear Anal., 11 (2022), 598-619.
\bibitem{Carvalho2017} Carvalho, M.L.M., da Silva, E.D.,  Goulart, C.: Quasilinear elliptic problems with concave-convex nonlinearities. Commun. Contemp. Math., 19 (2017), 1-25.
\bibitem{Silva2023} da Silva, E.D., Oliveira, J., Goulart, C.: Fractional $p$-Laplacian elliptic problems with sign changing nonlinearities via the nonlinear Rayleigh quotient. J. Math. Anal. Appl., 526 (2023).
\bibitem{Fan2014} Fan, H.N.: Multiple positive solutions for semi-linear elliptic systems with sign-changing weight. J. Math. Anal. Appl., 409 (2014), 399-408.
\bibitem{Fan2015} Fan, H.N., Liu, X.C.: Multiple positive solutions for semi-linear elliptic systems involving sign-changing weight. Math. Method Appl. Sci., 38 (2015), 1342-1351.
\bibitem{Fan2013} Fan, H.N., Liu, X.C.: Multiple positive solutions for degenerate elliptic equations with critical cone Sobolev exponents on singular manifolds. Electron. J. Differ. Eq., 181 (2013), 1-22.
\bibitem{Hamdani2021} Hamdani, M.K., Chung, N.T., Repov\v{s}, D.D.: New class of sixth-order nonhomogeneous $p(x)$-Kirchhff problems with sign-changing weight functions. Adv. Nonlinear Anal., 10 (2021), 1117-1131.
\bibitem{Huang2013} Huang, J.C., Chen, C.S., Xiu, Z.H.: Existence and multiplicity results for a $p$-Kirchhoff equation with a concave-convex term. Appl. Math. Lett., 26 (2013), 1070-1075.
\bibitem{Jeanjean2022} Jeanjean, L., R\u{a}dulescu, V.D.: Nonhomogeneous quasilinear elliptic problems: linear and sublinear cases. J. Anal. Math., 146 (2022), 327-350.
\bibitem{Lin2012} Lin, H.L.: Multiple positive solutions for semilinear elliptic systems. J. Math. Anal. Appl., 391 (2012), 107-118.
\bibitem{Li2013} Li, Q., Yang, Z.D.: Multiple positive solutions for quasilinear elliptic systems. Electron. J. Differ. Equations, 2013 (2013), 1-14.
\bibitem{Li2014} Li, Q., Yang, Z.D.: Multiple positive solutions for quasilinear elliptic systems with critical exponent and sign-changing weight. Comput. Math. Appl., 67 (2014), 1848-1863.
\bibitem{Landau1984} Landau, L.D., Lifshitz, E.M., Pitaevskii, L.P.: Electrodynamics of Continuous
Media. New York: Pergamon Press, 1984.
\bibitem{Mawhin1989} Mawhin, J., Willem, M.: Critical Point Theorem and Hamiltonian System. 1st edition, SpringerVerlag, New York, 1989.
\bibitem{Mihalache1989} Mihalache, D., Bertolotti, M., Sibilia, C.: Nonlinear wave propagation in planar structures.
Prog. Opt., 27 (1989), 227-313.
\bibitem{Pomponio2021} Pomponio, A., Watanabe, T.: Ground state solutions for quasilinear scalar field equations arising in nonlinear optics. NoDEA-Nonlinear Diff., 26 (2021).
\bibitem{Qi2023} Qi, W.T., Zhang, X.Y.: Multiplicity of solutions for a nonhomogeneous quasilinear elliptic
equation with concave-convex nonlinearitiesar. Xiv:2311.00283.
\bibitem{Rabinowitz1986} Rabinowitz, P.H.: Minimax Methods in Critical Point Theory with Applications to Differential Equations. 1st edition, American Mathmatical Society, the United States of America, 1986.
\bibitem{Sahu1991} Sahu, A., Priyadarshi, A.: Semilinear elliptic equation involving the $p$-Laplacian on the Sierpi\'{n}ski gasket. Complex Var. Ellipic, 64 (2019), 112-125.
\bibitem{Stuart2001} Stuart, C.A., Zhou, H.S.: Existence of guided cylindrical TM-modes in a homogeneous self-focusing dielectric. Ann. I. H. Poincar\'e-An., 18 (2001), 69-96.
\bibitem{Stuart2010} Stuart, C.A.,  Zhou, H.S.: Existence of guided cylindrical TM-modes in an inhomogeneous self-focusing dielectric. Math. Mod. Meth. Appl. S., 20 (2010), 1681-1719.
\bibitem{Stuart2011} Stuart, C.A.: Two positive solutions of a quasilinear elliptic Dirichlet problem. Milan J. Math., 79 (2011), 327-341.
\bibitem{Stuart1996} Stuart, C.A.: Cylindrical TM-modes in a homogeneous self-focusing dielectric.  Math. Mod. Meth. Appl. S., 6 (1996), 977-1008.
\bibitem{Stuart1997} Stuart, C.A.: Magnetic field wave equations for TM-modes in nonlinear optical waveguides.  Reaction Diffusion Systems, Editors Caristi and Mitidieri, Marcel Dekker, New York, 1997.
\bibitem{Svelto1974} Svelto, O.: Self-focusing, self-trapping and self-phase modulation in Laser beams.
Prog. Opt., 12 (1974), 1-51.
\bibitem{Wu2008}  Wu, T.F.: The Nehari manifold for a semilinear elliptic system involving sign-changing weight function. Nonlinear Anal-Theor., 68 (2008), 1733-1745.
\bibitem{Wu2006} Wu, T.F.: On semilinear elliptic equations involving concave-convex nonlinearities and sign-changing weight function. J. Math. Anal. Appl., 318 (2006), 253-270.
\bibitem{Wu2009}  Wu, T.F.: Multiplicity results for a semilinear elliptic equation involving sign-changing weight function. Rocky Mountain J. Math., 39 (2009), 995-1012.
\bibitem{Wu2010} Wu, T.F.: Multiple positive solutions for a class of concave-convex elliptic problems in $\mathbb{R}^{N}$ involving sign changing weight. J. Funct. Anal., (2010) 99-131.
\bibitem{Xiang2015} Xiang, M.Q., Zhang, B.L., Ferrara, M.: Multiplicity results for the nonhomogeneous fractional
    $p$-Kirchhoff equations with concave-convex nonlinearities. P. Roy. Soc. A-Math. Phy., 2015.
\bibitem{Xiu2016} Xiu, Z., Zhao, J., Chen, J., Li, S.: Multiple solutions on a $p$-biharmonic equation with nonlocal term. Bound value probl., 2016.
\bibitem{Zuo2019} Zuo, J., An, T., Ye, G., Qiao, Z.: Nonhomogeneous fractional $p$-Kirchhoff problems involving a critical nonlinearity. Electron. J. Qual. Theo., 41 (2019), 1-15.












\end{thebibliography}
\end{document}